\documentclass[a4paper,12pt]{amsart}
\usepackage[utf8]{inputenc}
\usepackage{amsmath,amssymb,amsthm}
\usepackage[top=2.8cm,bottom=2.5cm,left=2.6cm,right=2.6cm]{geometry}
\usepackage[numbers,compress, sort]{natbib}
\usepackage[colorlinks=true,
            linkcolor=blue,
            citecolor=blue,
            urlcolor=blue,
            filecolor=blue,
            menucolor=blue,
            runcolor=blue]{hyperref}
\usepackage{lmodern}
\usepackage{mathrsfs}
\usepackage[normalem]{ulem}
\usepackage{tikz-cd}
\usepackage{enumitem}
\usepackage{diagbox}
\usepackage{array}
\usepackage{booktabs}
\usepackage{multirow}
\usepackage{float}
\usepackage[utf8]{inputenc}
\usepackage{amsmath, amssymb}
\usepackage{hyperref}

\usepackage{parskip}
\makeatletter
\newtheoremstyle{mystyle}
{0.7em}               
{0.5em}               
{\itshape}            
{}                       
{\bfseries}              
{.}                      
{0.5em}                  
{}                       
\makeatother

\theoremstyle{mystyle}

\newtheorem{Theorem}{Theorem}[section]
\newtheorem*{Theorem*}{Theorem}
\newtheorem{Lemma}[Theorem]{Lemma}
\newtheorem{Corollary}[Theorem]{Corollary}
\newtheorem{Proposition}[Theorem]{Proposition}

\newtheorem{Remark}[Theorem]{Remark}

\theoremstyle{definition}
\newtheorem{Example}[Theorem]{Example}
\theoremstyle{plain}
\newenvironment{Proof}{\par\noindent%
\textit{\textbf{Proof.}}
}{\hfill $\square$\par}

\allowdisplaybreaks[4]

\makeatletter
\@namedef{subjclassname@2010}{{\rm2020} Mathematics Subject Classification}
\makeatother

\makeatletter

 \numberwithin{equation}{section}
\allowdisplaybreaks

\makeatother
\title[Isoparametric hypersurfaces in $\mathbb{S}^{n}\times \mathbb{R}^{m}$ and $\mathbb{H}^{n}\times \mathbb{R}^{m}$]{\textbf{Isoparametric hypersurfaces in $\mathbb{S}^{\lowercase{n}}\times \mathbb{R}^{\lowercase{m}}$ and $\mathbb{H}^{\lowercase{n}}\times \mathbb{R}^{\lowercase{m}}$}}
\author[H. X. Tan]{Huixin Tan\textsuperscript{1}}
\address{\textsuperscript{1}School of Mathematical Sciences, Laboratory of Mathematics and Complex Systems, Beijing Normal University, Beijing, 100875, P. R. China}
\email{hxtan@mail.bnu.edu.cn}
\author[Y. Q. Xie]{Yuquan Xie\textsuperscript{2}}
\address{\textsuperscript{2}School of Mathematics, Hangzhou Normal University, Hangzhou 311121, P. R. China}
\email{yuqxie@hznu.edu.cn}
\author[W. J. Yan]{Wenjiao Yan\textsuperscript{3,$\ast$}}
\address{\textsuperscript{3,$\ast$}School of Mathematical Sciences, Laboratory of Mathematics and Complex Systems, Beijing Normal University, Beijing, 100875, P. R. China}
\email{wjyan@bnu.edu.cn}

\thanks{$^{\ast}$ the corresponding author}

\subjclass[2010] { 53C42, 53B25, 53C40.}
\keywords{Isoparametric hypersurface, Homogeneous hypersurface, Product space}

\begin{document}

\begin{abstract}
    We first show that every isoparametric hypersurface in $\mathbb{S}^{n}\times \mathbb{R}^{m}$ or $\mathbb{H}^{n}\times \mathbb{R}^{m}$ possesses a constant angle function with respect to the canonical product structure. Exploiting this rigidity, we achieve a complete classification of isoparametric and homogeneous hypersurfaces in these product spaces. Furthermore, we prove that an isoparametric hypersurface in $\mathbb{S}^{n}\times \mathbb{R}^{m}$ or $\mathbb{H}^{n}\times \mathbb{R}^{m}$ also has constant principal curvatures.

\end{abstract}

\maketitle

\section{Introduction}\label{sec:intro}
A smooth non-constant function $F : M \to \mathbb{R}$ on a Riemannian manifold $M$ is called \emph{transnormal} if there exists a smooth function $b:\mathbb{R}\rightarrow \mathbb{R}$ such that $\|\nabla F\|^2 = b(F)$, where $\nabla F$ denotes the gradient of $F$. If, in addition, there exists a continuous function \( a : \mathbb{R} \to \mathbb{R} \) such that the Laplacian satisfies $\Delta F = a(F)$, then $F$ is said to be \emph{isoparametric} (cf. \cite{Wang-1987}). The regular level sets $\Sigma=F^{-1}(t)$
are correspondingly referred to as \emph{transnormal} or \emph{isoparametric hypersurfaces}, respectively. As observed by Élie Cartan, the transnormal condition implies that the level hypersurfaces are parallel, while the isoparametric condition further guarantees that these parallel hypersurfaces have constant mean curvatures. Moreover, in real space forms, Cartan proved that a hypersurface is isoparametric if and only if it has constant principal curvatures. 

The classification of isoparametric hypersurfaces in the Euclidean space $\mathbb{R}^{n}$ and hyperbolic space $\mathbb{H}^{n}$ was completed by Cartan \cite{Cartan-1938} and Segre \cite{Segre-1938} as early as in 1938. By contrast, the $\mathbb{S}^n$ case remained a subtle and long-standing problem—indeed, S. T. Yau listed it as Problem 34 in ``Open Problems in Geometry" \cite{Schoen-Yau-2010}.
After decades of contributions from numerous mathematicians \cite{Munzner-1980,Munzner-1981,Cartan-1939-a,Cartan-1939-b,Abresch-1983,Miyaoka-2013-Annals,Miyaoka-2016-Annals,Ozeki-Takeuchi-1975,Ozeki-Takeuchi-1976,Ferus-Karcher-Munzner-1981,Solomon-1992,Grove-Halperin-1978,Tang-1991,Fang-1999,Stolz-1999,Cecil-Chi-Jensen-2007,Immervoll-2008,Chi-2011,Chi-2013}, a complete classification on the unit sphere $\mathbb{S}^n$ was finally achieved in 2020 \cite{Chi-2020-JDG}.
Furthermore,
the theory of isoparametric hypersurfaces serves as a powerful tool for investigating diverse geometric problems \cite{Ge-Xie-2010, Tang-Yan-2013, Tang-Yan-2015, Tang-Yan-2014, Tang-Yan-2022, Tang-Yan-2023, Tang-Yan-2025}.
A natural continuation of this classical theme is to study the classification of isoparametric hypersurfaces in the Riemannian product of two real space forms, $M_{c_{1}}^{n} \times M_{c_{2}}^{m}$ ($c_1, c_2\in \{1, 0, -1\}$).

In order to classify isoparametric and homogeneous hypersurfaces in the product manifold $\mathbb{S}^{2} \times \mathbb{S}^{2}$, Urbano \cite{Urbano-2019} introduced in 2019 a natural product structure $P$ on the tangent bundle of $\mathbb{S}^{2} \times \mathbb{S}^{2}$, together with an associated angle function $C$ defined on an oriented hypersurface $\Sigma$. These constructions, in fact, extend verbatim to any product of two real space forms
$M_{c_{1}}^{n} \times M_{c_{2}}^{m}$ ($c_1, c_2\in \{1, 0, -1\}$). 
Concretely, if a tangent vector decomposes as  $(v_{1}, v_{2})$ according to the product splitting, the structure $P$ is defined by:
\begin{equation*}
	\begin{aligned}
		P: \mathfrak{X} (M_{c_{1}}^{n} \times M_{c_{2}}^{m}) ~&\longrightarrow~ \mathfrak{X} (M_{c_{1}}^{n} \times M_{c_{2}}^{m}) \\
		(v_{1}, v_{2}) ~&\longmapsto~ (v_{1}, -v_{2}),
	\end{aligned}
\end{equation*}
where $\mathfrak{X}(\cdot)$ denotes the set of all smooth tangent vector fields on the manifold.
With respect to the product metric, this tensor field satisfies $P^2=Id$, and is parallel. Let $\Sigma\subset M_{c_{1}}^{n} \times M_{c_{2}}^{m}$ be an orientable hypersurface with the unit normal vector field $N$. With respect to the product metric, the associated \emph{angle function} $C$ is defined by
\begin{equation*}
	\begin{aligned}
		C: \Sigma~ &\longrightarrow~ [-1, 1] \\
		x ~&\longmapsto ~\langle PN(x), N(x) \rangle,
	\end{aligned}
\end{equation*}
which measures the projection of the normal vector onto the $\pm 1$-eigenspaces of $P$.
The extreme values $C=\pm 1$ correspond to normals entirely contained in one factor, while $|C|<1$ indicates a genuine tilt between the two factors.

Recent works establish the following rigidity for isoparametric hypersurfaces:

\vspace{3mm}
\noindent
\textbf{Theorem} (\cite{Urbano-2019, Gao-Ma-Yao-2024-a, Gao-Ma-Yao-2024-b, de-Lima-Pipoli-2024})
	\textit{In each of the spaces
	$\mathbb{S}^{2} \times \mathbb{S}^{2}$,
	$\mathbb{S}^{2} \times \mathbb{R}^{2}$,
	$\mathbb{S}^{2} \times \mathbb{H}^{2}$,
	$\mathbb{H}^{2} \times \mathbb{H}^{2}$,
	$\mathbb{H}^{2} \times \mathbb{R}^{2}$,
	$\mathbb{S}^{n} \times \mathbb{R}$,
	and $\mathbb{H}^{n} \times \mathbb{R}$ \((\)with $n \ge 2$\()\),
	all isoparametric hypersurfaces have the constant angle.}
\vspace{3mm}

In this paper, we extend these results to higher-dimension Euclidean factors:

\begin{Theorem}\label{thm-main}
  Let $\Sigma$ be a connected isoparametric hypersurface in $\mathbb{S}^{n} \times \mathbb{R}^{m}$ or $\mathbb{H}^{n} \times \mathbb{R}^{m}$.
  Then the associated angle function $C$ is constant along $\Sigma$.
\end{Theorem}

\begin{Remark}
In a forthcoming paper, we shall classify isoparametric hypersurfaces in the remaining product spaces $\mathbb{S}^{n} \times \mathbb{S}^{m}$, $\mathbb{S}^{n} \times \mathbb{H}^{m}$, and $\mathbb{H}^{n} \times \mathbb{H}^{m}$, and establish a corresponding constant-angle property. These investigations will provide a foundation for broader geometric problems, such as those in \cite{Ding-Ge-Li-2025, Li-Wang-2024, Xia-Xiao-Zhong-2025}.
\end{Remark}

\begin{Remark}\label{remark-main-thm-proof}
	The case $n \ge 2$ in Theorem \ref{thm-main} is proved in Section \ref{sec-proof-thm-main}.
	The argument used there does not apply when $n = 1$; nevertheless, Example \ref{ex-S1Rm} together with Theorem \ref{thm-class-isop-MnRm}-(iii) yields a direct verification that $\Sigma$ has the constant angle function in the $n=1$ case.
\end{Remark}

Urbano \cite{Urbano-2019} obtained a complete classification of isoparametric hypersurfaces in \(\mathbb{S}^{2}\times\mathbb{S}^{2}\) by constructing an efficient global frame adapted to the natural complex structures on \(\mathbb{S}^{2}\). Several subsequent works followed his strategy to treat other product models. In 2018 Julio–Batalla \cite{Julio-Batalla-2018} classified isoparametric hypersurfaces with constant principal curvatures in \(\mathbb{S}^{2}\times\mathbb{R}^{2}\); later, dos Santos–dos Santos \cite{dos-Santos-dos-Santos-2023} treated the case \(M_{c_1}^2\times M_{c_2}^2\) with \(c_1\neq c_2\). Gao–Ma–Yao \cite{Gao-Ma-Yao-2024-b} removed the assumption of constant principal curvatures in \cite{dos-Santos-dos-Santos-2023} and completed the classification; in a related work \cite{Gao-Ma-Yao-2024-a} they developed refined geometric tools to treat \(\mathbb{H}^{2}\times\mathbb{H}^{2}\). All these approaches crucially exploit the fact that every two-dimensional real space form carries a natural complex structure, and therefore their arguments do not generalize to higher dimensions (for instance, among all spheres only \(\mathbb{S}^2\) and \(\mathbb{S}^6\) admit almost complex structures, while whether \(\mathbb{S}^6\) carries a complex structure remains the well-known Hopf problem \cite{Tang-Yan-2022,Tang-Yan-2025}).

It is also noteworthy that Ge-Radeschi \cite{Ge-Radeschi-2015} obtained a foliated diffeomorphism classification of codimension one singular Riemannian foliations (e.g. isoparametric foliation) on all closed simply connected $4$-manifolds (including $\mathbb{S}^{2}\times\mathbb{S}^{2}$).
In addition, Qian-Tang \cite{Qian-Tang-2016} provided an isoparametric hypersurface in $\mathbb{S}^{n}\times\mathbb{S}^{n}$ and computed its curvature properties as well as the spectrum of the Laplace–Beltrami operator.
More recently, Cui \cite{Cui-2025} provided further examples of isoparametric hypersurfaces by restricting certain isoparametric functions on $\mathbb{S}^{2n+1}$ to the product $\mathbb{S}^{n} \times \mathbb{S}^{n}$.

From another perspective, building upon the local classification of constant angle hypersurfaces in~\cite{de-Lima-Roitman-2021}, de Lima and Pipoli~\cite{de-Lima-Pipoli-2024} obtained a complete classification of isoparametric hypersurfaces in 
\(\mathbb{S}^n\times\mathbb{R}\) and \(\mathbb{H}^n\times\mathbb{R}\). 
They proved the following result:

\vspace{2mm}
\noindent
\textbf{Theorem} (\cite{de-Lima-Pipoli-2024}) \textit{Isoparametric hypersurfaces in $M_{c}^{n}\times \mathbb{R}~(c=\pm 1)$ are precisely one of the following:
\begin{enumerate}[label=\textup{(\roman*)}, itemsep=-2pt, topsep=0pt]
	\item horizontal slice $M_{c}^{n}\times \{t_{0}\}$;
	\item a vertical cylinder over a complete isoparametric hypersurface in $M_{c}^{n}$;
	\item a parabolic bowl in $\mathbb{H}^{n}\times\mathbb{R}$.
\end{enumerate}}
\vspace{2mm}

The classification above is based on the concept of \((M_s,\phi)\)-graphs. 
However, a direct extension of this construction to vector-valued functions 
produces submanifolds of higher codimension rather than hypersurfaces, 
and thus does not apply when the Euclidean factor has dimension \(m>1\).

We adopt a different approach. Inspired by Miyaoka \cite{Miyaoka-2013} and through a focused analysis along the special principal direction $V=PN-CN$, we establish the following classification of isoparametric hypersurfaces in \(M_c^n\times\mathbb{R}^m\) \((c=\pm1)\). (Note that when $n=1$, only the case 
\(\mathbb{S}^1\times\mathbb{R}^m\) needs to be considered, as \(\mathbb{H}^1\) does not exist.)
\vspace{1mm}

\begin{Theorem}\label{thm-class-isop-MnRm}
	Let \(\Sigma\) be a connected complete isoparametric hypersurface in \(M_c^n\times\mathbb{R}^m\) \((c=\pm1,\ m\ge 2)\), i.e., in \(\mathbb{S}^n\times\mathbb{R}^m\) or \(\mathbb{H}^n\times\mathbb{R}^m\). Up to ambient isometry, \(\Sigma\) is one of the following:
	
	\begin{enumerate}[label=\textup{(\roman*)}, itemsep=-1pt, topsep=1pt]
		\item \(K_{1}\times\mathbb{R}^m\), where \(K_{1}\) is an isoparametric hypersurface in \(M_c^n\), while for \(n=1\) this reduces to \(\{p\}\times\mathbb{R}^m\), \(p\in\mathbb{S}^1\);
		\item \(M_c^n\times K_{2}\), where \(K_{2}\) is an isoparametric hypersurface in \(\mathbb{R}^m\);
		\item \(\Phi(\mathbb{R}^m)\subset\mathbb{S}^1\times\mathbb{R}^m\), where \(\Phi\colon\mathbb{R}^m\to\mathbb{S}^1\times\mathbb{R}^m\) is the immersion defined by
		\[
		x\mapsto\big(\cos\langle x,x_0\rangle,\ \sin\langle x,x_0\rangle,\ x\big),
		\]
		with \(\langle\cdot,\cdot\rangle\) denoting the standard inner product on \(\mathbb{R}^m\) and \(x_0\in\mathbb{R}^m\setminus\{0\}\) fixed;
		\item \(\Psi(\mathbb{R}^{n+m-1})\subset\mathbb{H}^n\times\mathbb{R}^m\), where \(\Psi\colon\mathbb{R}^{n+m-1}\to\mathbb{H}^n\times\mathbb{R}^m\) is given by
	$$(t,x,y)\mapsto\big(p(t,x),\,q(t,y)\big),$$
		with
            \begin{align*}
        	p(t,x) &= \cosh \big(t\sqrt{\varepsilon}\big) \gamma_{1}(x) + \sinh \big(t\sqrt{\varepsilon} \big) N_{\gamma_{1}}(x), \\
        	q(t,y) &= \gamma_{2}(y) + t\sqrt{1-\varepsilon} N_{\gamma_{2}},
           \end{align*}
		where \(\gamma_1(x)\) is a horosphere in \(\mathbb{H}^n\) with the unit normal \(N_{\gamma_1}\), \(\gamma_2(y)\) is an affine hyperplane in \(\mathbb{R}^m\) with the constant unit normal \(N_{\gamma_2}\), and \(\varepsilon\in (0, 1)\) is a constant.
	\end{enumerate}
\end{Theorem}
 \vspace{1mm}

\begin{Remark}
In a forthcoming paper, we shall generalize this classification to the remaining product types \(\mathbb{S}^n\times\mathbb{S}^m\), \(\mathbb{S}^n\times\mathbb{H}^m\) and \(\mathbb{H}^n\times\mathbb{H}^m\).
\end{Remark}

 \vspace{1mm}

As mentioned earlier, in real space forms, isoparametric hypersurfaces coincide with hypersurfaces having constant principal curvatures. However, these two notions are no longer equivalent in general Riemannian manifolds. For example, Rodr{\'i}guez-V{\'a}zquez \cite{Rodriguez-Vazquez-2019} constructed non-isoparametric hypersurfaces with constant principal curvatures in the torus \(\mathbb{T}^{n}\) (\(n \ge 3\)), while Ge-Tang-Yan \cite{Ge-Tang-Yan-2015} exhibited isoparametric hypersurfaces in \(\mathbb{C}P^n\) whose principal curvatures are not constant.

Combining Theorems \ref{thm-main}, \ref{thm-class-isop-MnRm}, with Example \ref{ex-S1Rm},  Example \ref{ex-HnRm}, we immediately obtain the following characterization.
\vspace{1mm}

\begin{Corollary}\label{coro-main-equiv}
	Let \(\Sigma\) be a connected complete isoparametric hypersurface in \(\mathbb{S}^{n}\times \mathbb{R}^{m}\) or \(\mathbb{H}^{n}\times \mathbb{R}^{m}\).  
	Then \(\Sigma\) has constant principal curvatures.
\end{Corollary}
 \vspace{1mm}

Furthermore, by combining Theorem \ref{thm-class-isop-MnRm}, Corollary \ref{coro-main-equiv} with Example \ref{ex-S1Rm},  Example \ref{ex-HnRm}, we obtain a classification of homogeneous hypersurfaces in \(\mathbb{S}^{n}\times\mathbb{R}^{m}\) and \(\mathbb{H}^{n}\times\mathbb{R}^{m}\). This result generalizes that of~\cite{de-Lima-Pipoli-2024}, which corresponds to the case $m=1$.
\vspace{1mm}

\begin{Corollary}\label{coro-class-homo-MnRm}
	Let \(\Sigma\) be a homogeneous hypersurface in \(M_{c}^{n}\times\mathbb{R}^{m}\) \( (c=\pm1\), \(m\ge 2 )\), i.e., in \(\mathbb{S}^{n}\times\mathbb{R}^{m}\) or \(\mathbb{H}^{n}\times\mathbb{R}^{m}\).  
	Up to ambient isometry, \(\Sigma\) is one of the following:
	
	\begin{enumerate}[label=\textup{(\roman*)}, itemsep=-2pt, topsep=0pt]
		\item \(K_{1}\times\mathbb{R}^{m}\), where \(K_{1}\) is a homogeneous hypersurface in \(M_{c}^{n}\), while in the case \(n=1\), this reduces to \(\{p\}\times\mathbb{R}^{m}\) with \(p\in\mathbb{S}^{1}\);
		\item \(M_{c}^{n}\times K_{2}\), where \(K_{2}\) is a homogeneous hypersurface in \(\mathbb{R}^{m}\);
		\item The hypersurface described in Theorem \ref{thm-class-isop-MnRm}-(iii);
		\item The hypersurface described in Theorem \ref{thm-class-isop-MnRm}-(iv).
	\end{enumerate}
\end{Corollary}

\vspace{1mm}

The paper is organized as follows.  
In Section \ref{sec-proof-thm-class}, we prove Theorem \ref{thm-class-isop-MnRm} and verify the homogeneity of the hypersurfaces listed therein.  
Due to its length and technical nature, the proof of Theorem \ref{thm-main} is presented separately in Section \ref{sec-proof-thm-main}.

\vspace{3mm}



\section{Preliminaries}\label{sec:preliminaries}

Let $\Sigma$ be an orientable hypersurface in the product manifold $M_{c_{1}}^{n} \times M_{c_{2}}^{m}$ with the global unit normal vector field $N$.  
For any vector field $X \in \mathfrak{X}(M_{c_{1}}^{n} \times M_{c_{2}}^{m})$, we denote by $X^{h}$ its \emph{horizontal component} tangent to $M_{c_{1}}^{n}$ and by $X^{v}$ its \emph{vertical component} tangent to $M_{c_{2}}^{m}$. 
Let $A$ be the shape operator of $\Sigma$ associated with $N$,
and $H$ the mean curvature of $\Sigma$. The natural projection maps are given by
\begin{align*}
	\pi_{1}:M_{c_{1}}^{n} \times M_{c_{2}}^{m}~
	&\longrightarrow ~M_{c_{1}}^{n},
	&
	\pi_{2}:M_{c_{1}}^{n} \times M_{c_{2}}^{m}~
	&\longrightarrow ~M_{c_{2}}^{m}, \\
	(x,y)~
	&\longmapsto ~x,
	&
	(x,y)~
	&\longmapsto~ y.
\end{align*}
For each $(x,y) \in \Sigma$, we define
\[
\Sigma_{y} = \pi_{1}\big(\pi_{2}^{-1}(y) \cap \Sigma\big)
\quad\text{and}\quad
\Sigma_{x} = \pi_{2}\big(\pi_{1}^{-1}(x) \cap \Sigma\big),
\]
which represent the projections of $\Sigma$ into the horizontal and vertical factors, respectively.

Decompose the unit normal vector as $N = (N^{h}, N^{v})$. The angle function 
$C$ is then accordingly given by
\begin{equation}\label{equ-C=PNN}
	C = \langle P N, N \rangle 
	= \|N^{h}\|^{2} - \|N^{v}\|^{2}
	= C_{1}^{2} - C_{2}^{2},
\end{equation}
where
\begin{equation*}
C_{1} = \|N^{h}\| = \sqrt{\frac{1+C}{2}},
\quad
C_{2} = \|N^{v}\| = \sqrt{\frac{1-C}{2}}.
\end{equation*}

Now we introduce a special tangent vector field $V$ on $\Sigma$, which will play an important role in the subsequent verification. It is defined by
\begin{equation}\label{equ-def-V}
	V = P N - C N=\left((1-C)N^h, -(1+C)N^v\right).
\end{equation}
It follows immediately that $\|V\|^{2} = 1 - C^{2}$.
Differentiating \eqref{equ-C=PNN} and using the fact that $P$ is parallel, we obtain, for any tangent vector field $X$ on $\Sigma$,
\begin{align*}
	X(C)
	&= \langle \nabla_{X}(P N), N \rangle
	+ \langle P N, \nabla_{X} N \rangle \\
	&= -2 \langle A X, V \rangle
	= -2 \langle X, A V \rangle.
\end{align*}
Hence, the gradient of $C$ is given by
\begin{equation}\label{equ-nablaC=-2AV}
	\nabla^{\Sigma} C = -2 A V.
\end{equation}

\medskip

In the product $M_{c}^{n} \times \mathbb{R}^{m} (c=\pm 1)$, the Riemannian curvature tensor $R_{c}$ of the product manifold $M_{c}^{n} \times \mathbb{R}^{m}$ is given by
\begin{equation}\label{equ-curvature_tensor}
	R_{c}(X,Y)Z 
	= c \big( \langle X^{h}, Z^{h} \rangle Y^{h}
	- \langle Y^{h}, Z^{h} \rangle X^{h} \big),
	\qquad
	\forall\, X,Y,Z \in \mathfrak{X}(M_{c}^{n} \times \mathbb{R}^{m}).
\end{equation}


\section{Classification of Isoparametric Hypersurfaces}\label{sec-proof-thm-class}

In this section, we aim to prove Theorem \ref{thm-class-isop-MnRm}.
For the fluency of expression, we begin by preparing two propositions to characterize the focal points and principal frames of transnormal hypersurfaces with a constant angle in general Riemannian product manifolds $M_{1}\times M_{2}$.
 \vspace{1mm}
 
\begin{Proposition}\label{prop-Nk0-const-tildeF-trans}
	Let $\Sigma$ be a connected complete transnormal hypersurface in the Riemannian product $M_{1}\times M_{2}$.
	If the angle function $C$ is constant with $-1<C<1$, then for any $(x_{0},y_{0})\in \Sigma$, the slices $\Sigma_{x_{0}}$ and $\Sigma_{y_{0}}$ are transnormal hypersurfaces in $M_{2}$ and $M_{1}$, respectively.
\vspace{-1mm}	

	Moreover, if $(x, y) \in M_{1} \times M_{2}$ is a focal point of $\Sigma$, then $x$ is a focal point in $M_{1}$ and $y$ is a focal point in $M_{2}$. Conversely, if $x$ is a focal point in $M_{1}$ or $y$ is a focal point in $M_{2}$, then $(x, y)$ is a focal point in $M_{1} \times M_{2}$.
\end{Proposition}
 \vspace{1mm}

\begin{Proof}
Without loss of generality, let $\Sigma=F^{-1}(t)$ be a regular level set of a transnormal function $F: M_{1} \times M_{2} \rightarrow \mathbb{R}$ satisfying $\|\nabla F\|^{2} = b(F)$.
Denote by $\nabla^{h}$ and $\nabla^{v}$ the gradients on $M_{1}$ and $M_{2}$, respectively.
	
For fixed points $x_0\in M_1$ and $y_0\in M_2$, define
\begin{align*}
	\begin{aligned}
		F_{x_{0}}:M_{2}
		&\longrightarrow \mathbb{R}, \\
		y
		&\longmapsto F(x_{0},y),
	\end{aligned}
	\qquad
	\begin{aligned}
		F_{y_{0}}:M_{1}
		&\longrightarrow \mathbb{R}, \\
		x
		&\longmapsto F(x,y_{0}).
	\end{aligned}
\end{align*}
Then $\Sigma_{x_{0}} = F_{x_{0}}^{-1}(t)$ and $\Sigma_{y_{0}} = F_{y_{0}}^{-1}(t)$.
A straightforward computation yields
	  \begin{align*}
		\lVert\nabla^{v} F_{x_{0}}(y)\rVert^2  
		&= \lVert \nabla^{v} F(x_{0},y)\rVert^2 \\
		&= \frac{1-C}{2} \|\nabla F(x_{0},y)\rVert^{2}
		= \frac{1-C}{2} b\big(F(x_{0},y)\big),
	\end{align*}
and the corresponding relation for $F_{y_{0}}$ is analogous. Moreover, by (\ref{equ-def-V})
   \begin{align*}
	\exp_{(x_{0},y_{0})}\frac{2}{1-C}t(0,N^{v})
	&= \exp_{(x_{0},y_{0})}\frac{1}{1-C}t\big((1-C)N-V\big) \\
	&= \exp_{\exp_{(x_{0},y_{0})}\left(-\frac{1}{1-C}tV\right)} tN,
\end{align*}
and similarly,
\begin{align*}
	\exp_{(x_{0},y_{0})}\frac{2}{1+C}t(N^{h},0)
	&= \exp_{(x_{0},y_{0})}\frac{1}{1+C} t\big((1+C)N+V\big) \\
	&= \exp_{\exp_{(x_{0},y_{0})}\frac{1}{1+C}tV}tN.
\end{align*}

Since $-1<C<1$, we have $\Vert V \Vert^{2} = 1-C^{2} \neq 0$. Hence $V$ is a nonvanishing tangent vector field on the complete hypersurface $\Sigma$. Therefore, $\exp_{(x_{0},y_{0})} tV$ defines a diffeomorphism for each $t$, 
and the differential of $\exp_{y_{0}} t N^{v}$ (resp., $\exp_{x_{0}} t N^{h}$) has the same rank as that of $\exp_{(x_{0},y_{0})} tN$. The desired conclusion follows.
	
\end{Proof}
 \vspace{1mm}
 
\begin{Remark}
	When $C \equiv 1$, we have $N = (N^{h}, 0)$.
	Thus, for any $(x_{0}, y_{0}) \in \Sigma$, $\Sigma_{y_{0}} = M_{1}$ and $\Sigma_{x_{0}}$ is a transnormal hypersurface in $M_{2}$.
	Consequently, $(x,y) \in M_{1} \times M_{2}$ is a focal point if and only if $y \in M_{2}$ is a focal point. The case $C \equiv -1$ is analogous.
\end{Remark}
 \vspace{1mm}

\begin{Proposition}\label{prop-AXiXj-AYiYj}
	Let $\Sigma$ be a connected complete transnormal hypersurface with the constant angle function $C$ in a Riemannian product $M_{1}^{n}\times M_{2}^{m}$, and set $V=PN-CN$. 
	If $-1<C<1$, then at any point $(x,y)\in \Sigma$, there exists a local orthonormal frame
	\begin{align*}
		\left\{ \frac{1}{\sqrt{1-C^{2}}}V, ~(X_{1}, 0), \ldots, (X_{n-1}, 0), ~(0, Y_{n}), \ldots, (0, Y_{n+m-2}) \right\}
	\end{align*}
with respect to which the shape operator $A$ of $\Sigma$ satisfies
\begin{equation}\label{equ-Sigma-Principla-Frame}
\begin{cases}
	\qquad \qquad AV&=0, \\
	\langle A(X_{i}, 0), (X_{j}, 0) \rangle 
	&= \lambda_{i} \delta_{ij},
	\quad i,j=1,\ldots,n-1, \\
	\langle A(0, Y_{\alpha}), (0, Y_{\beta}) \rangle 
	&= \lambda_{\alpha} \delta_{\alpha \beta},
	\quad \alpha,\beta=n,\ldots,n+m-2.
\end{cases}
\end{equation}
Here, $\lambda_{i}/C_{1}$ $(i=1,\ldots,n-1)$ are the principal curvatures of $\Sigma_{y}$ in $M_{1}^{n}$, and $\lambda_{\alpha}/C_{2}$ $(\alpha=n,\ldots,n+m-2)$ are the principal curvatures of $\Sigma_{x}$ in $M_{2}^{m}$. Moreover, the mean curvature of $\Sigma$ is given by $H=\sum_{i=1}^{n+m-2}\lambda_{i}$.
\end{Proposition}
 \vspace{1mm}
 
\begin{Proof}
Denote by $\nabla$, $\nabla^{h}$, and $\nabla^{v}$ the Levi-Civita connections on $M_{1}^{n}\times M_{2}^{m}$, $M_{1}^{n}$, and $M_{2}^{m}$, respectively.
Since the hypersurface $\Sigma$ has a constant angle function $C$, for any $X\in \mathfrak{X}(\Sigma_{y})$ and $Y\in \mathfrak{X}(\Sigma_{x})$, we have
\begin{align*}
	\langle\nabla_{X}^{h}N^{h},N^{h}\rangle_{M_{1}}
	&=\frac{1}{2}X\langle N^{h},N^{h}\rangle_{M_{1}}
	=0, \\
	\langle\nabla^{v}_{Y}N^{v},N^{v}\rangle_{M_{2}}
	&=\frac{1}{2}Y\langle N^{v},N^{v}\rangle_{M_{2}}
	=0.
\end{align*}
Let $A^{h}$ denote the shape operator of $\Sigma_{y}\subset M_{1}^{n}$, such that for the unit normal vector field $\eta=N^{h}/C_{1}$ and any tangent vector $X$, we have $A^{h}_{\eta}X=-\nabla_{X}^{h}\eta$. Similarly, let $A^{v}$ denote the shape operator of $\Sigma_{x}\subset M_{2}^{m}$. Then
\begin{equation*}
	\begin{aligned}
		C_{1}A^{h}X
		&=-\nabla_{X}^{h}N^{h}+\langle\nabla_{X}^{h}N^{h},N^{h}\rangle_{M_{1}} N^{h}
		=-\nabla_{X}^{h}N^{h}, \\
		C_{2}A^{v}Y
		&=-\nabla_{Y}^{v}N^{v}+\langle\nabla_{Y}^{v}N^{v},N^{v}\rangle_{M_{2}} N^{v}
		=-\nabla_{Y}^{v}N^{v}.
	\end{aligned}
\end{equation*}
	
At each point $(x,y)\in \Sigma$, let $\{(X_{i},0)\}_{i=1}^{n-1}$ be eigenvectors of $A^{h}$ corresponding to eigenvalues $\lambda_{i}/C_{1}$, and $\{(0,Y_{\alpha})\}_{\alpha=n}^{n+m-2}$ be eigenvectors of $A^{v}$ corresponding to eigenvalues $\lambda_{\alpha}/C_{2}$, respectively. Then we have
\begin{align*}
    \langle A(X_{i},0),(X_{j},0)\rangle
    &=-\langle\nabla_{(X_{i},0)}(N^{h},N^{v}),(X_{j},0)\rangle
    =-\langle\nabla_{X_{i}}^{h}N^{h},X_{j} \rangle_{M_{1}^{n}} \\
    &=\langle C_{1}A^{h}X_{i},X_{j}\rangle_{M_{1}^{n}} 
    =\lambda_{i}\delta_{ij}
\end{align*}
for any $i,j=1,\ldots,n-1$ and similarly,
\begin{align*}
    \langle A(0,Y_{\alpha}),(0,Y_{\beta})\rangle
    &=-\langle \nabla_{(0,Y_{\alpha})}(N^{h},N^{v}),(0,Y_{\beta})\rangle 
    =-\langle\nabla_{Y_{\alpha}}^{v}N^{v},Y_{\beta} \rangle_{M_{2}^{m}} \\
    &=\langle C_{2}A^{v}Y_{\alpha},Y_{\beta}\rangle_{M_{2}^{m}} 
    =\lambda_{\alpha}\delta_{\alpha \beta}
\end{align*}
for any $\alpha,\beta=n,\ldots,n+m-2$.     
Therefore, equation \eqref{equ-Sigma-Principla-Frame} follows.

\end{Proof}

\vspace{5mm}
Now, we proceed to complete the proof of Theorem \ref{thm-class-isop-MnRm}.
\vspace{3mm}

\noindent
\textit{Proof of Theorem \ref{thm-class-isop-MnRm}.$\quad$}
According to Theorem \ref{thm-main}, the isoparametric hypersurface $\Sigma$ possesses a constant angle function $C$. The cases $C=1$ and $C=-1$ correspond to (i) and (ii), respectively; hence we assume $-1 < C < 1$ in the sequel.

We first consider the case $n=1$, i.e., $\mathbb{S}^1 \times \mathbb{R}^m$, which leads to parts (i)–(iii) of the classification.
Recall the following result from \cite{Ge-Tang-2013}.
 \vspace{1mm}
 
\begin{Lemma}[\cite{Ge-Tang-2013}]\label{lem-Ge-Tang-2013-prop2.9}
	Let $\pi:E\to B$ be a Riemannian submersion with minimal fibers.
	Given any (properly) isoparametric function $f$ on $B$,
	then $F:=f\circ \pi$ is a (properly) isoparametric function on $E$.
\end{Lemma}
 \vspace{1mm}

The universal cover $\pi:\mathbb{R}\to\mathbb{S}^1$, $\pi(x)=\mathrm{e}^{\sqrt{-1}x}$, has discrete (hence minimal) fibers, and the induced covering map
    \begin{align*}
	\widetilde{\pi}: \mathbb{R}^{m+1} &\longrightarrow \mathbb{S}^{1} \times \mathbb{R}^{m} \\ (x,y)&\longmapsto (\pi(x),y)
\end{align*}
is a Riemannian submersion with minimal fibers. By Lemma \ref{lem-Ge-Tang-2013-prop2.9}, it  suffices to find an isoparametric function $F$ on $\mathbb{R}^{m+1}$ satisfying $F = f \circ \widetilde{\pi}$, where $f$ is an isoparametric function on $\mathbb{S}^{1} \times \mathbb{R}^{m}$. 
Notice that the periodicity of $\widetilde{\pi}$ implies that
$F(x+2k\pi,y)=F(x,y)$ for all $x\in \mathbb{R}$ and $k\in\mathbb{Z}$.

If the foliation determined by $F$ admits a focal manifold $\Sigma_0$, the classification in $\mathbb{R}^{m+1}$ implies that $\Sigma_0$ is either a single point or an affine subspace of dimension at most $m-1$. Moreover, for any $(x,y)\in\Sigma_0$, the entire line $\mathbb{R}\times\{y\}\subset\Sigma_0$; otherwise $\Sigma_0$ would decompose into disjoint union of lower-dimensional affine subspaces, which does not occur in the classification.    Consequently, $F(x + 2k\pi, y) = F(x, y)$ for any $y \in \mathbb{R}^{m}$, and the identity $F(x, y) = F(x', y)$ holds for all $x, x' \in \mathbb{R}$, thereby proving Theorem \ref{thm-class-isop-MnRm}-(ii).

If $F$ admits no focal manifold, the classification in $\mathbb{R}^{m+1}$ implies that its regular level sets must be hyperplanes. The periodicity condition allows one to choose
$F(x,y) = \sin\left(x - \kappa \langle y, y_{0} \rangle\right)$, where $y_{0}$ is a unit vector in $\mathbb{R}^{m}$. 
When $\kappa=0$, a connected component of the regular level set of $F$ corresponds to Theorem \ref{thm-class-isop-MnRm}-(i); when $\kappa \neq 0$, each connected component of a regular level set of $F$ can be parameterized as in Theorem \ref{thm-class-isop-MnRm}-(iii).

Next, consider $n\ge2$, which leads to parts (i), (ii), and (iv). The constancy of $C$ implies that $C_1$ and $C_2$ are also constant. Denote by $F$ the isoparametric function on $M_c^n\times\mathbb{R}^m$ associated with $\Sigma$.
\medskip

\noindent
\textbf{Case 1: $\mathbb{S}^{n}\times\mathbb{R}^{m}$ $(n\geq 2)$. }
For any $(x,y) \in \Sigma$, Proposition \ref{prop-Nk0-const-tildeF-trans} shows that $\Sigma_{x}$ and $\Sigma_{y}$ are regular level sets of transnormal functions on $\mathbb{R}^{m}$ and $\mathbb{S}^{n}$, respectively, and hence are isoparametric by \cite[Theorem 1.5-(1)]{Miyaoka-2013}.
However, isoparametric hypersurfaces in $\mathbb{S}^n$ have focal points that occur infinitely often along each normal geodesic. By applying
\begin{align*} 
	\exp_{(x,y)}tN = \left( x \cos C_{1} t + \frac{\sin C_{1} t}{C_{1}} N^{h}, \, y + t N^{v} \right),
\end{align*} 
it follows that $\Sigma_x$ would have infinitely many focal points in $\mathbb{R}^m$, contradicting their classification.
Therefore, no isoparametric hypersurfaces with $-1<C<1$ exist in $\mathbb{S}^n\times\mathbb{R}^m$.
\medskip

\noindent
\textbf{Case 2: $\mathbb{H}^{n}\times\mathbb{R}^{m}$.}
By \cite[Theorem 1.1]{Miyaoka-2013}, the possible topological types of $\mathbb{H}^{n}\times\mathbb{R}^{m}$ are as follows:
\begin{enumerate}[label=(\roman*), itemsep=-1pt, topsep=1pt]
    \item if the transnormal system has no focal submanifold: an $\mathbb{R}$-bundle or $\mathbb{S}^{1}$-bundle over a hypersurface $\Sigma$,
    \item if there is one focal submanifold: either a vector bundle over the unique focal submanifold $\widetilde{\Sigma}$ or a DDBD structure,
    \item if there are two focal submanifolds: a DDBD structure,
\end{enumerate}
where the DDBD (Double Disc Bundle Decomposition) structure means that the ambient manifold is constructed by glueing two disc bundles over two submanifolds along the boundaries.

The $\mathbb{S}^1$-bundle case is excluded since $\exp_{(x,y)}tN\ne(x,y)$ for any $t\ne0$.
	
The DDBD structure is also impossible. If $\mathbb{H}^{n}\times\mathbb{R}^{m}$ admitted a DDBD structure, then for any point $(x,y)\in\Sigma$, the normal geodesic would intersect the focal manifold infinitely many times, yielding infinitely many focal points along it. By Proposition \ref{prop-Nk0-const-tildeF-trans}, this implies that $\Sigma_{x}$ also has infinitely many focal points along the normal geodesic in $\mathbb{R}^{m}$, contradicting the known focal structure of isoparametric hypersurfaces in $\mathbb{R}^{m}$.

As for the remaining two cases, we first show that $\Sigma_y$ is an isoparametric hypersurface in $\mathbb{H}^{n}$, since the isoparametricity of $\Sigma_{x}$ in $\mathbb{R}^{m}$ follows from a similar discussion as in Case 1. Furthermore, we will see that $\Sigma_{x}$ is isometric to $\Sigma_{x'}$ and $\Sigma_{y}$ is isometric to $\Sigma_{y'}$ for any $(x,y),(x',y')\in\Sigma$.

In case that $\mathbb{H}^{n} \times \mathbb{R}^{m}$ is a vector bundle over its unique focal submanifold $\widetilde{\Sigma}$, let $\Sigma\subset\mathbb{H}^{n}\times\mathbb{R}^{m}$ be the tube of constant radius $t$ around $\widetilde{\Sigma}$. 
For any $(x,y)\in\Sigma$, we have
\begin{align*} 
\exp_{\exp_{(x,y)}(-\tfrac{t}{1-C}V)}tN
=& \exp_{(x,y)}\tfrac{2t}{1-C}(0,N^{v}), \\
\exp_{\exp_{(x,y)}(\tfrac{t}{1+C}V)}tN
=&\exp_{(x,y)}\tfrac{2t}{1+C}(N^{h},0).
\end{align*}
Hence $\Sigma_{x}\subset\mathbb{R}^{m}$ lies at distance $\tfrac{t}{C_{2}}$ from its focal submanifold along the unit normal $N^{v}/C_{2}$, 
and $\Sigma_{y}\subset\mathbb{H}^{n}$ lies at distance $\tfrac{t}{C_{1}}$ along $N^{h}/C_{1}$. 
Moreover, $\Sigma_{x}$ has constant mean curvature $H_{\Sigma_{x}}=\ell C_{2}/t$ for some integer $\ell \in \{1, \ldots, m-1\}$.
By Proposition \ref{prop-AXiXj-AYiYj}, since $\Sigma$ has a constant angle function $C\in (-1, 1)$,
\[
H_{\Sigma}(x,y) = C_{1} H_{\Sigma_{y}}(x) + C_{2} H_{\Sigma_{x}}(y),
\]
and since $H_{\Sigma}$, $H_{\Sigma_{x}}$, $C_1$ and $C_2$ are constant, so is $H_{\Sigma_{y}}$.
Let $\Sigma_{t}$ denote the parallel hypersurface at distance $t$ from $\Sigma$, and $\Sigma_{y,t}$ the parallel hypersurface at distance $t$ from $\Sigma_{y}\subset \mathbb{H}^{n}$. Noting that
\begin{align*} 
\exp_{(x,y)} \frac{2}{1+C} t (N^{h},0) = \exp_{\exp_{(x,y)} \frac{1}{1+C} t V} t N,
\end{align*}
we obtain
\begin{align*} 
H_{\Sigma_{C_{1} t}}(x,y) = C_{1} H_{\Sigma_{y,t}}\Big(\exp_{x} t \frac{N^{h}}{C_{1}}\Big) + C_{2} H_{\Sigma_{\exp_{x} t \frac{N^{h}}{C_{1}}}}(y).
\end{align*}
Since $H_{\Sigma_{C_{1} t}}$ and $H_{\Sigma_{\exp_{x} t \frac{N^{h}}{C_{1}}}}$
are both constant, so is $H_{\Sigma_{y,t}}$.
Hence $\Sigma_{y}$ is an isoparametric hypersurface in $\mathbb{H}^{n}$. Moreover, since an isoparametric hypersurface with a single focal submanifold in $\mathbb{H}^{n}$ or $\mathbb{R}^{m}$ is uniquely determined (up to isometry) by its distance to the focal submanifold, 
it follows that $\Sigma_{x}$ is isometric to $\Sigma_{x'}$ and $\Sigma_{y}$ to $\Sigma_{y'}$ for any $(x,y),(x',y')\in\Sigma$.

In case that $\mathbb{H}^{n}\times \mathbb{R}^{m}$ is an $\mathbb{R}$-bundle over the hypersurface $\Sigma$, $\Sigma\subset\mathbb{R}^{m}$ has no focal points. It follows that $\Sigma_{x}$ also has no focal points, and thus is a hyperplane with vanishing mean curvatures. Then an analogous discussion as the previous case shows that $\Sigma_{y}\subset\mathbb{H}^{n}$ is isoparametric.
Moreover, for any $(x,y)\in\Sigma$, $\Sigma_{y}\subset\mathbb{H}^{n}$ is one of the following:
\begin{enumerate}[label=(\roman*), itemsep=0pt, topsep=1pt]
	\item a totally geodesic hyperplane ($\lambda_{i}=0$),
	\item an equidistant hypersurface ($0<|\lambda_{i}|<C_{1}$), or
	\item a horosphere ($\lambda_{i}=\pm C_{1}$),
\end{enumerate}
while $\Sigma_{x}\subset\mathbb{R}^{m}$ is a hyperplane ($\lambda_{\alpha}=0$). In these cases, principal curvatures at $(x,y),(x',y')\in\Sigma$ are the same. 
Since isoparametric hypersurfaces in $\mathbb{H}^{n}$ and $\mathbb{R}^{m}$ without focal submanifolds are uniquely determined by their principal curvatures, 
we again conclude that $\Sigma_{x}$ and $\Sigma_{y}$ are pairwise isometric.

\medskip

Next, we consider the flow along $V$:
\[
\exp_{(x,y)}tV
= \big(\exp_{x}\!\big((1-C)tN^{h}\big),\, \exp_{y}\!\big(-(1+C)tN^{v}\big)\big).
\]
From Proposition~\ref{prop-AXiXj-AYiYj}, we have
\begin{align*}
	A(X_{i},0) &= \lambda^{H}_{i}(X_{i},0) + \sigma_{i\alpha}(0,Y_{\alpha}), \\ 
	A(0,Y_{\alpha}) &= \sigma_{\alpha i}(X_{i},0) + \lambda^{R}_{\alpha}(0,Y_{\alpha}), 
\end{align*}
where $(\sigma_{i\alpha})$ is an $(n-1)\times(m-1)$ matrix.
Let $A_{t}$ denote the shape operator at $\exp_{(x,y)}tV$.  
Since $\Sigma_{y}$ and $\Sigma_{x}$ are isometric to $\Sigma_{\exp_{y}(-(1+C)tN^{v})}$ and $\Sigma_{\exp_{x}((1-C)tN^{h})}$, respectively, 
we may assume
\begin{align*}
	A_{t}(X_{i},0) &= p_{ki}\lambda^{H}_{k}p_{kj}(X_{j},0)
	+ p_{ki}\sigma_{k\alpha}q_{\alpha\beta}(0,Y_{\beta}), \\ 
    \quad
	A_{t}(0,Y_{\alpha}) &= q_{\gamma\alpha}\sigma_{i\gamma}p_{ij}(X_{j},0)
	+ q_{\gamma\alpha}\lambda^{R}_{\gamma}q_{\gamma\beta}(0,Y_{\beta}), 
\end{align*}
where $(p_{ij})$ and $(q_{\alpha\beta})$ are orthogonal matrices of orders $n-1$ and $m-1$.
Taking $(X,Y)=(X_{i},0)$ and $(X,Y)=(0,Y_{\alpha})$ in the equation $A_{t}\big((f_{t})_{*}(X,Y)\big)=-\nabla_{(f_{t})_{*}(X,Y)}N$ reveals the following formulas
\begin{align}
\label{equ-HnRm1-Xi1}
&R_{1}(t) (p_{ij})^{\mathrm{T}} \Lambda_{H} (p_{ij}) +(1+C) t (\sigma_{ij})(q_{\alpha \beta})^{\mathrm{T}} (\sigma_{\alpha i}) (p_{ij}) = -\frac{1}{1-C} R'_{1}(t), \\
\label{equ-HnRm1-Yi1}
&(I +(1+C)t \Lambda_{R}) (q_{\alpha \beta})^{\mathrm{T}} \Lambda_{R} (q_{\alpha \beta}) +R_{2}(t) (\sigma_{\alpha i}) (p_{ij})^{\mathrm{T}} (\sigma_{i\alpha}) (q_{\alpha \beta}) = \Lambda_{R},
\end{align}
where
\begin{align*}
R_{1}(t) &= \cosh C_{1}(1-C) t I -\frac{1}{C_{1}}\sinh C_{1}(1-C) t \Lambda_{H}, \\
R_{2}(t) &= -\frac{1}{C_{1}} \sinh C_{1}(1-C) t,
\end{align*}
and $I$ denotes the identity matrix,
$\Lambda_{H}$ is the diagonal matrix composed of $\lambda^{H}_{i}$,
and $\Lambda_{R}$ is the diagonal matrix composed of $\lambda^{R}_{\alpha}$.
Multiplying both sides of equation \eqref{equ-HnRm1-Xi1} (resp. \eqref{equ-HnRm1-Yi1}) from the right by the matrix $(p_{ij})^{\mathrm{T}}$ (resp. $(q_{\alpha \beta})^{\mathrm{T}}$) and utilizing the orthogonality of $(p_{ij})$ (resp. $(q_{\alpha \beta})$),
we obtain
\begin{align*}
    R_{1}(t)
    (p_{ij})^{\mathrm{T}}
    \Lambda_{H}
    +(1+C) t(\sigma_{i\alpha})
    (q_{\alpha \beta})^{\mathrm{T}}(\sigma_{\alpha i})
    =&
    -\frac{1}{1-C} R'_{1}(t)
    (p_{ij})^{\mathrm{T}}, \\
    \left(I+(1+C)t\Lambda_{R} \right)(q_{\alpha \beta})^{\mathrm{T}} \Lambda_{R}
    +R_{2}(t) (\sigma_{\alpha i}) (p_{ij})^{\mathrm{T}} (\sigma_{i\alpha})
    =&\Lambda_{R} (q_{\alpha \beta})^{\mathrm{T}}.
\end{align*}
Differentiating both sides of the above equations with respect to $t$ at $t=0$ and focusing on the diagonal entries, we obtain
\[
\left(\lambda^{H}_{i}\right)^{2} = C_{1}^{2} + \frac{1+C}{1-C}\sum_{\alpha}\sigma_{i\alpha}^{2}, \qquad \left(\lambda^{R}_{\alpha}\right)^{2} = \frac{1-C}{1+C}\sum_{i}\sigma_{\alpha i}^{2}.
\]
Summing over the indices $i$ and $\alpha$, respectively, leads to
\begin{align}
\sum_{i}(\lambda^{H}_{i})^{2} &= (n-1)C_{1}^{2} + \frac{(1+C)^{2}}{(1-C)^{2}}\sum_{\alpha}(\lambda^{R}_{\alpha})^{2}. \label{equ-HnRm1-PH-QR}
\end{align}


If $\Sigma_{y}$ and $\Sigma_{x}$ focalize simultaneously, we can view $\Sigma$ as a tube of radius $s$ around the focal submanifold, giving
\begin{align*}
	\lambda^{H}_{1} = \cdots = \lambda^{H}_{k} &= C_{1} \coth \frac{s}{C_{1}}, &
	\lambda^{H}_{k+1} = \cdots = \lambda^{H}_{n-1} &= C_{1} \tanh \frac{s}{C_{1}},\\
	\lambda^{R}_{1} = \cdots = \lambda^{R}_{\ell} &= \frac{C_{2}^{2}}{s}, &
	\lambda^{R}_{\ell+1} = \cdots = \lambda^{R}_{m-1} &= 0,
\end{align*}
for some $k\in \{1,\dots,n-1\}$ and $\ell \in \{1,\dots,m-1\}$. However, equation \eqref{equ-HnRm1-PH-QR} contradicts the above equations for all $s$.

If $\Sigma$ has no focal points, the classification implies that $\Sigma_{x}$ is a hyperplane in $\mathbb{R}^{m}$, i.e., $\lambda^{R}_{i}=0$. Substituting this into \eqref{equ-HnRm1-PH-QR}, we obtain that $\lambda^{H}_{1} = \cdots = \lambda^{H}_{n-1} = \pm C_{1}$. Therefore, the only remaining case is that $\Sigma_{x}$ is a hyperplane in $\mathbb{R}^{m}$ and $\Sigma_{y}$ is a horosphere in $\mathbb{H}^{n}$ for any $(x,y)\in \Sigma$, with $\sigma_{i\alpha}=0$, leading directly to the expression in Theorem \ref{thm-class-isop-MnRm}-(iv).

\hfill$\Box$
\vspace{3mm}

\begin{Example}\label{ex-S1Rm}
	From the proof of Theorem~\ref{thm-class-isop-MnRm}, the hypersurfaces described in case~(iii) arise as connected components of the level sets of
	\[
	F:\mathbb{S}^{1}\times\mathbb{R}^{m} \to \mathbb{R}, \qquad
	F(e^{\sqrt{-1}x},y) = \sin\big(x - \kappa \langle y, y_{0} \rangle\big),
	\]
	where $y_{0}$ is a fixed unit vector in $\mathbb{R}^{m}$ and $\kappa\in\mathbb{R}$.
	
	A direct computation yields
	\[
	\nabla F = \big(\cos(x - \kappa \langle y, y_{0} \rangle),\, -\kappa \cos(x - \kappa \langle y, y_{0} \rangle)y_{0}\big),
	\]
	and hence
	\[
	\|\nabla F\|^{2} = (1+\kappa^{2})(1-F^{2}), \qquad
	\Delta F = -(1+\kappa^{2})F.
	\]
	Thus, $F$ is an isoparametric function on $\mathbb{S}^{1}\times\mathbb{R}^{m}$.
	
	All values of $F$ except $\pm1$ are regular. Moreover, $F^{-1}(\pm1)$ are also connected isoparametric hypersurfaces parameterized as in Theorem~\ref{thm-class-isop-MnRm}-(iii). For each $t\in(-1,1)$, $F^{-1}(t)$ consists of two connected components.
	
	\medskip
	
	Consider $\Sigma = F^{-1}(t)$ for $t \in [-1,1]$. Its unit normal vector is
	\[
	N = \mathrm{sgn}\big(\cos(x - \kappa \langle y, y_{0} \rangle)\big)
	\frac{1}{\sqrt{1+\kappa^{2}}}(1, -\kappa y_{0}),
	\]
	and the angle function is \( C = \frac{1-\kappa^{2}}{1+\kappa^{2}}. \)
	A straightforward computation shows that the Hessian of $F$ satisfies $\nabla^{2}F|_{\Sigma}=0$, implying that $\Sigma$ is totally geodesic in $\mathbb{S}^{1}\times\mathbb{R}^{m}$.
	
	Indeed,
	\[
	\nabla^{2}F
	=\sin\big(x-\kappa\langle y,y_{0} \rangle\big)
	\begin{pmatrix}
		-1 & \kappa y_{0}^{\mathrm{T}}\\
		\kappa y_{0} & -\kappa^{2}y_{0}y_{0}^{\mathrm{T}}
	\end{pmatrix}.
	\]
	Choose an orthonormal frame $\{v_{1},\ldots,v_{m}\}$ on $\Sigma$, where $v_{i}=(0,Y_{i})$ for $i=1,\ldots,m-1$ and
	\[
	v_{m}=\frac{V}{\|V\|}
	=\mathrm{sgn}\big(\cos(x-\kappa\langle y,y_{0}\rangle)\big)
	\Big(\tfrac{|\kappa|}{\sqrt{1+\kappa^{2}}},\,\tfrac{\mathrm{sgn}(\kappa)}{\sqrt{1+\kappa^{2}}}y_{0}\Big).
	\]
	For $i,j\le m-1$, since $Y_{j}\perp y_{0}$ and $\|y_0\|=1$, one verifies
	\[
	\nabla^{2}F|_{\Sigma}(v_{i},v_{j})=0,\qquad
	\nabla^{2}F|_{\Sigma}(v_{i},v_{m})=0,
	\]
	and hence $\nabla^{2}F|_{\Sigma}=0$.
	
	\medskip
	
	We now show that each connected component of a level set of $F$ is homogeneous in $\mathbb{S}^{1}\times\mathbb{R}^{m}$.  
	When $\kappa=0$, $F=\sin x$, corresponding to case~(i) of Theorem~\ref{thm-class-isop-MnRm}. Hence we assume $\kappa\ne0$.
	
	Let $\Sigma_{0}$ be a connected component of $\Sigma=F^{-1}(t)$ for $t\ne\pm1$.
    Denote by $\mathrm{Isom}_{0}(\mathbb{S}^{1}\times\mathbb{R}^{m})$ the identity component of the isometry group of $\mathbb{S}^{1}\times\mathbb{R}^{m}$, which is isomorphic to $SO(2) \times (\mathbb{R}^{m} \rtimes SO(m))$, where $SO(2)$ is the special orthogonal group of degree $2$ and $\mathbb{R}^{m} \rtimes SO(m)$ is the special Euclidean group in $m$ dimensions.    
	Represent $y\in\mathbb{R}^{m}$ by $\begin{pmatrix}y\\1\end{pmatrix}$ and consider the subgroup
	\[
	K=\langle K_{1},K_{2}\rangle
	\subset\mathrm{Isom}_{0}(\mathbb{S}^{1}\times\mathbb{R}^{m}),
	\]
	where
	\[
	K_{1}=\Bigg\{
	\begin{pmatrix}
		1 & 0 & 0\\
		0 & B & b\\
		0 & 0 & 1
	\end{pmatrix}
	\!\!\Bigg|
	B\in SO(m),\,B^{\mathrm{T}}y_{0}=y_{0},\,\langle b,y_{0}\rangle=0
	\Bigg\},
	\]
	and
	\[
	K_{2}=\Bigg\{
	\begin{pmatrix}
		e^{\sqrt{-1}\theta} & 0 & 0\\
		0 & I_{m} & \tfrac{\theta}{\kappa}y_{0}\\
		0 & 0 & 1
	\end{pmatrix}
	\!\!\Bigg|
	\theta\in\mathbb{R}
	\Bigg\}.
	\]
	Since $K_{1}$ and $K_{2}$ commute and $K_{1}\cap K_{2}=\{Id\}$, we have $K\cong K_{1}\times K_{2}$, i.e.,
	\[
	K=\Bigg\{
	\begin{pmatrix}
		e^{\sqrt{-1}\theta} & 0 & 0\\
		0 & B & b+\tfrac{\theta}{\kappa}y_{0}\\
		0 & 0 & 1
	\end{pmatrix}
	\!\!\Bigg|
	B^{\mathrm{T}}y_{0}=y_{0},\,
	\langle b,y_{0}\rangle=0,\,
	B\in SO(m),\,\theta\in\mathbb{R}
	\Bigg\}.
	\]
	Define
    \begin{align*}
        \phi: \mathrm{Isom}_{0}(\mathbb{S}^{1}\times\mathbb{R}^{m})
        &\longrightarrow SO(2)\times\mathbb{R}^{2}, \\
        \left(e^{\sqrt{-1}\theta},\begin{pmatrix}
            B & b \\
            0 & 1
        \end{pmatrix}
        \right)
        &\longmapsto \Big(e^{\sqrt{-1}\theta},\langle B^{\mathrm{T}}y_{0},y_{0}\rangle,\langle b,y_{0}\rangle\Big).
    \end{align*}
	Evidently, $\phi$ is continuous. Since $K=\phi^{-1}(D)$ with
	$D=\{(e^{\sqrt{-1}\theta},1,\theta/\kappa)\mid\theta\in\mathbb{R}\}$ closed in $SO(2)\times\mathbb{R}^{2}$, $K$ is a closed subgroup of $\mathrm{Isom}(\mathbb{S}^{1}\times\mathbb{R}^{m})$.
	
	\medskip
	
	Finally, $\Sigma_{0}$ is an orbit of $K$.  
	For $(x,y),(x',y')\in\Sigma_{0}$, we have $x'-x=\kappa\langle y'-y,y_{0}\rangle$. Moreover, since $y$ and $y' + \langle y - y', y_{0} \rangle y_{0}$ lie in the same hyperplane perpendicular to $y_{0}$, we can choose $B \in SO(m)$ and $b\in\mathbb{R}^{m}$ such that $B^{\mathrm{T}} y_{0} = y_{0}$ and $By+b+\frac{x'-x}{\kappa}y_{0}=y'$. Then
	\[
	\begin{pmatrix}
		e^{\sqrt{-1}(x'-x)} & 0 & 0\\
		0 & B & b+\tfrac{x'-x}{\kappa}y_{0}\\
		0 & 0 & 1
	\end{pmatrix}\!\!\in K
	\]
	maps $(x,y)$ to $(x',y')$. Thus, $K$ acts transitively on $\Sigma_{0}$ and preserves it, proving that $\Sigma_{0}$ is a homogeneous hypersurface in $\mathbb{S}^{1}\times\mathbb{R}^{m}$.  
	The same argument applies when $t=\pm1$.
	
	\hfill$\Box$
\end{Example}

\vspace{3mm}

\begin{Example}\label{ex-HnRm}
	The isoparametric function corresponding to case (iv) of Theorem~\ref{thm-class-isop-MnRm} is
	\[
	F:\mathbb{H}^{n}\times\mathbb{R}^{m}\to\mathbb{R}, \qquad
	F(x,y)=\langle x,u\rangle_{\mathbb{L}}\exp\big(a\langle y-y_{0},v_{0}\rangle\big),
	\]
	where $u=(u_{0},\ldots,u_{n})$ is a nonzero lightlike vector in the Lorentz space $\mathbb{L}^{n+1}$ with $u_{0}>0$, $\langle\cdot,\cdot\rangle_{\mathbb{L}}$ denotes the Lorentz inner product, $v_{0}$ is a fixed unit vector in $\mathbb{R}^{m}$, $y_{0}\in\mathbb{R}^{m}$, and $a\in\mathbb{R}$.
	
	A direct computation gives
	\[
	\nabla F=\big((u+\langle x,u\rangle_{\mathbb{L}}x)\exp(a\langle y-y_{0},v_{0}\rangle),\; a v_{0}\langle x,u\rangle_{\mathbb{L}}\exp(a\langle y-y_{0},v_{0}\rangle)\big),
	\]
	and hence
	\[
	\|\nabla F\|^{2}=(1+a^{2})F^{2}, \qquad \Delta F=(n+a^{2})F.
	\]
	Thus, $F$ is an isoparametric function, and all its level sets are regular.
	
	For fixed $y\in\mathbb{R}^m$, the equation $F(x,y)=t$ gives $\langle x,u\rangle_{\mathbb{L}}=t\,\exp(-a\langle y-y_{0},v_{0}\rangle)$, representing a horosphere in $\mathbb{H}^{n}$ centered at the lightlike vector $u$. For fixed $x\in\mathbb{H}^n$, one obtains $\langle y-y_{0},v_{0}\rangle=\frac{1}{a}\ln\frac{t}{\langle x,u\rangle_{\mathbb{L}}}$, defining an affine hyperplane in $\mathbb{R}^{m}$ through $y_{0}$ with the unit normal $v_{0}$.
	
	Let $\Sigma=F^{-1}(t)$ for $t\in(-\infty,0)$. Its unit normal and angle function are
	\[
	N=\frac{1}{\sqrt{1+a^{2}}}\Big(\frac{u+\langle x,u\rangle_{\mathbb{L}}x}{\langle x,u\rangle_{\mathbb{L}}},\; a v_{0}\Big),
	\qquad
	C=\frac{1-a^{2}}{1+a^{2}}.
	\]
	The cases $a=0$ and $|a|\to\infty$ correspond to Theorem~\ref{thm-class-isop-MnRm}-(i) and -(ii), respectively; thus we focus on $a\neq0$.
	
  For tangent vectors $X=(X^{h},X^{v})$ and $Y=(Y^{h},Y^{v})$ of $\mathbb{H}^n\times\mathbb{R}^m$, the Hessian of $F$ is
	\begin{align*}
    	\nabla^{2} F(X,Y)
    	=& \langle X^{h}, Y^{h} \rangle F + a^2 \langle X^{v}, v_0 \rangle \langle Y^{v}, v_0 \rangle F \\
    	&+ a \exp\big(a\langle y-y_{0},v_{0}\rangle \big) \big( \langle X^{h}, u^{\top} \rangle \langle Y^{v}, v_0 \rangle + \langle Y^{h}, u^{\top} \rangle \langle X^{v}, v_0 \rangle \big),
    \end{align*}
    where $u^{\top}=u+\langle u,x\rangle_{\mathbb{L}}x$ is the projection of $u$ onto $T_{x}\mathbb{H}^{n}$.
	
	Choose an orthonormal frame $\{(X_{1},0),\dots,(X_{n-1},0),(0,Y_{1}),\dots,(0,Y_{m-1}),V/\|V\|\}$ on $\Sigma$, where $V=PN-CN$ and
	\[
	\frac{V}{\|V\|}
	=\frac{1}{\sqrt{1+a^{2}}}
	\Big(|a|\frac{u+\langle x,u\rangle_{\mathbb{L}}x}{\langle x,u\rangle_{\mathbb{L}}},\; -\mathrm{sgn}(a)v_{0}\Big).
	\]
	Under this frame,
	\[
	\nabla^{2}F|_{\Sigma}=
	\begin{pmatrix}
		t\,I_{n-1}&0&0\\
		0&0&0\\
		0&0&0
	\end{pmatrix}.
	\]
Since the second fundamental form $\mathrm{II}=-\frac{1}{\|\nabla F\|}\nabla^{2}F|_{\Sigma}$, the principal distributions are
\begin{align*}
	\mathcal{V}_{1} &=\mathrm{span}\,\{(X_{i},0) \mid i=1,\ldots,n-1\}, \\
	\mathcal{V}_{2} &=\mathrm{span}\,\{(0,Y_{j}) \mid j=1,\ldots, m-1 \}, \\
	\mathcal{V}_{3} &=\mathrm{span}\,\left\{V\right\}.
\end{align*}
	with the corresponding principal curvatures and multiplicities:
	\begin{table}[H]\centering
		\begin{tabular}{ccc}
			\toprule
			Distribution & Principal curvature & Multiplicity\\
			\midrule
			$\mathcal{V}_{1}$ & $\dfrac{1}{\sqrt{1+a^{2}}}$ & $n-1$\\
			$\mathcal{V}_{2}$ & $0$ & $m-1$\\
			$\mathcal{V}_{3}$ & $0$ & $1$\\
			\bottomrule
		\end{tabular}
	\end{table}
	Hence the mean curvature of $\Sigma$ is $H=\dfrac{n-1}{\sqrt{1+a^{2}}}$. For any principal directions $X,Y$, the sectional curvature is
	\[
	K_{\Sigma}(X,Y)
	=K_{\mathbb{H}^{n}\times\mathbb{R}^{m}}(X,Y)
	+\frac{\mathrm{II}(X,X)\mathrm{II}(Y,Y)-\mathrm{II}(X,Y)^{2}}
	{\langle X,X\rangle\langle Y,Y\rangle-\langle X,Y\rangle^{2}},
	\]
	yielding the following table of sectional curvatures:
	\begin{table}[H]\centering
		\begin{tabular}{c|ccc}
			\diagbox{$X\in$}{$Y\in$} & $\mathcal{V}_{1}$ & $\mathcal{V}_{2}$ & $\mathcal{V}_{3}$\\
			\hline
			$\mathcal{V}_{1}$ & $-\dfrac{a^{2}}{1+a^{2}}$ & $0$ & $-\dfrac{a^{2}}{1+a^{2}}$\\
			$\mathcal{V}_{2}$ & $0$ & $0$ & $0$\\
			$\mathcal{V}_{3}$ & $-\dfrac{a^{2}}{1+a^{2}}$ & $0$ & $-$\\
		\end{tabular}
	\end{table}
	Thus the Ricci and scalar curvatures are
	\[
	\mathrm{Ric}_{\Sigma}X=
	\begin{cases}
		-\dfrac{(n-1)a^{2}}{1+a^{2}}, & X\in\mathcal{V}_{1}\cup\mathcal{V}_{3},\\[3mm]
		\qquad 0, & X\in\mathcal{V}_{2},
	\end{cases}
    \qquad R=-\frac{n(n-1)a^{2}}{1+a^{2}}.
	\]
	Clearly, $\Sigma$ is not Einstein when $a\neq0$.
	
	\medskip
	To show $\Sigma$ is homogeneous, define the subgroup
	\[
	G=\langle G_{1},G_{2}\rangle
	\subset\mathrm{Isom}_{0}(\mathbb{H}^{n}\times\mathbb{R}^{m})
	\cong SO^{+}(1,n)\times(\mathbb{R}^{m}\rtimes SO(m)),
	\]
	where $SO^{+}(1,n)$ denotes the identity component of the Lorentz group and
	\begin{align*}
		G_{1}&=\Bigg\{
		\begin{pmatrix}B&0&0\\0&I_{m}&s v_{0}\\0&0&1\end{pmatrix}
		\!\Bigg|\;
		B\in SO^{+}(1,n),\;B^{\mathrm{T}}u=e^{-as}u,\;s\in\mathbb{R}
		\Bigg\},\\
		G_{2}&=\Bigg\{
		\begin{pmatrix}I_{n}&0&0\\0&\widetilde{B}&b\\0&0&1\end{pmatrix}
		\!\Bigg|\;
		\widetilde{B}\in SO(m),\;\widetilde{B}^{\mathrm{T}}v_{0}=v_{0},\;
		\langle b,v_{0}\rangle=0
		\Bigg\}.
	\end{align*}
An analogous discussion as in Example \ref{ex-S1Rm} shows that $G_{1}$ and $G_{2}$ commute, and thus $G\cong G_{1}\times G_{2}$, i.e.,
	\[
	G=\Bigg\{
	\begin{pmatrix}
		B&0&0\\
		0&\widetilde{B}&b+s v_{0}\\
		0&0&1
	\end{pmatrix}
	\!\!\Bigg|\!
	\begin{array}{l}
		B^{\mathrm{T}}u=e^{-as}u,\ B\in SO^{+}(1,n),\ s\in\mathbb{R},\\
		\widetilde{B}^{\mathrm{T}}v_{0}=v_{0},\ \langle b,v_{0}\rangle=0,\ \widetilde{B}\in SO(m)
	\end{array}
	\Bigg\}.
	\]
    Define
	\[
	\eta:\mathrm{Isom}_{0}(\mathbb{H}^{n}\times\mathbb{R}^{m})
	\to\mathbb{R}^{3},\quad
	\Big(B,\begin{pmatrix}\widetilde{B}&\widetilde{b}\\0&1\end{pmatrix}\Big)
	\mapsto
	\big(\langle B^{\mathrm{T}}u,u\rangle_{\mathbb{L}},\langle\widetilde{B}^{\mathrm{T}}v_{0},v_{0}\rangle,\langle \widetilde{b},v_{0}\rangle\big).
	\]
    Evidently, $\eta$ is continuous. Then $G=\eta^{-1}(D)$ with $D=\{(e^{-as},1,s)\mid s\in\mathbb{R}\}$ closed in $\mathbb{R}^{3}$; hence $G$ is a closed subgroup of $\mathrm{Isom}_{0}(\mathbb{H}^{n}\times\mathbb{R}^{m})$.
	
	Finally, for $(x,y),(x',y')\in\Sigma=F^{-1}(t)\subset \mathbb{H}^{n}\times\mathbb{R}^{m}$, 
	$\langle x',u\rangle_{\mathbb{L}}=\langle x,u\rangle_{\mathbb{L}} \exp\big(-a\langle y'-y,v_{0} \rangle\big) $. Then the transitivity of the isometric $SO^{+}(1,n)$-action on $\mathbb{H}^n$ yields the existence of $B_{0}\in SO^{+}(1,n)$ such that $B_{0}x=x'$, and thus 
	\[
	\langle B_{0}x,u\rangle_{\mathbb{L}}
    =\langle x, B_{0}^{\mathrm{T}}u\rangle_{\mathbb{L}}
    = \langle x,u\rangle_{\mathbb{L}} \exp\big(-a\langle y'-y,v_{0} \rangle\big),
	\]
	which implies $B_{0}^{\mathrm{T}}u= \exp\big(-a\langle y'-y,v_{0} \rangle\big) u$.
	Similarly, there exist $\widetilde{B}_{0}\in SO(m)$ and $b_{0}\in\mathbb{R}^{m}$ such that $\widetilde{B}_{0}y+b_{0}+\langle y'-y,v_{0}\rangle v_{0}=y'$, $\langle b_{0},v_{0}\rangle=0$, and thus $\widetilde{B}_{0} v_{0} = v_{0}$.  
	Then
	\[
	g=\begin{pmatrix}B_{0}&0&0\\0&\widetilde{B}_{0}&b_{0}+\langle y'-y,v_{0}\rangle v_{0}\\0&0&1\end{pmatrix}\in G
	\]
	maps $(x,y)$ to $(x',y')$. Thus $G$ acts transitively on $\Sigma$, and since $F$ is $G$-invariant, $G$ preserves $\Sigma$. Therefore, $\Sigma$ is a homogeneous hypersurface in $\mathbb{H}^{n}\times\mathbb{R}^{m}$.
	
	\hfill$\Box$
\end{Example}


\section{Proof of Theorem \ref{thm-main}}\label{sec-proof-thm-main}

As noted in Remark \ref{remark-main-thm-proof}, we restrict attention to the case $n \ge 2$.
Since the angle function is continuous, it suffices to show that it is locally constant.
Hence, we only consider the case $-1<C<1$.
For convenience, set $C_{1} = \tau$.
Clearly,
\[
    C_{2} = \sqrt{1 - \tau^{2}}, \qquad
    C = 2\tau^{2} - 1, \qquad
    0 < \tau < 1.
\]

We will choose an orthonormal frame along the parallel hypersurface of $\Sigma$ and compute the coefficient matrix of the Jacobi field with respect to this frame.
Then, by analyzing the linear system satisfied by the mean curvature of the parallel hypersurface and its derivatives, we derive a nontrivial algebraic equation in $\tau$, which in turn shows that the angle function $C$ must be constant.

\medskip

Let $N_{p}$ denote the unit normal vector of $\Sigma$ at $p \in \Sigma$, and define the normal exponential map $\Phi_{r}: \Sigma \to M_{c}^{n} \times \mathbb{R}^{m}$ by  
$\Phi_{r}(p) = \exp_{p}(r N_{p})$.  
Then there exists a sufficiently small $\delta > 0$ such that, for all $r \in (-\delta, \delta)$,  
the map $\Phi_{r}$ is well defined and $\Sigma_{r} = \Phi_{r}(\Sigma)$ is an embedded hypersurface in $M_{c}^{n} \times \mathbb{R}^{m}$ at distance $r$ from $\Sigma$.  
Fix $p \in \Sigma$, and let $\gamma_{p}(r)$, $r \in (-\delta, \delta)$, be the geodesic in $M_{c}^{n} \times \mathbb{R}^{m}$ satisfying $\gamma_{p}(0) = p$ and $\gamma_{p}'(0) = N_{p}$.  
The vector field $N(r) = \gamma_{p}'(r)$ along $\gamma_{p}$ is parallel, and hence remains normal to $\Sigma_{r}$ at $\gamma_{p}(r)$.

We now choose unit orthonormal vector fields $U_{1}(r), \ldots, U_{m}(r)$ parallel along $\gamma_{p}$ such that the horizontal components of $U_{i}(r)$ $(i=1,\ldots,m-1)$ vanish, i.e., $U^{h}_{i}(r)=0$, and  
\[
    U_{m}(r) = \left( -\frac{C_{2}}{C_{1}} N^{h}(r), \frac{C_{1}}{C_{2}} N^{v}(r) \right).
\]
Together with $N(r)$, we extend these to obtain a unit orthonormal parallel frame  
\[
    N(r),\; U_{1}(r),\; \ldots,\; U_{n+m-1}(r)
\]
along $\gamma_{p}$.  
By orthogonality, for $i = m+1, \ldots, n+m-1$, the vector fields $U_{i}(r)$ have vanishing vertical components.

For each $j = 1, \ldots, n+m-1$, let $\zeta_{j}(r)$ be the Jacobi field along $\gamma_{p}$ satisfying
\[
    \zeta_{j}(0) = U_{j}(0), \qquad
    \zeta'_{j}(0) = -A U_{j}(0),
\]
and
\begin{equation}\label{equ-Jacobi-equation}
    \zeta''_{j} + R_{c}(\gamma_{p}', \zeta_{j})\gamma_{p}' = 0,
\end{equation}
where the Riemann curvature tensor $R_{c}$ is defined in \eqref{equ-curvature_tensor}.  
To compute \eqref{equ-Jacobi-equation}, we decompose $\zeta_{j}(r)$ in the orthonormal frame $\{U_{i}(r)\}_{i=1}^{n+m-1}$ as
\[
    \zeta_{j}(r) = \sum_{i=1}^{n+m-1} b_{ij}(r) U_{i}(r),
\]
where $b_{ij}(r)$ are smooth functions on $(-\delta, \delta)$ for $j = 1, \ldots, n+m-1$.  
Meanwhile, the shape operator $A$ with respect to the orthonormal basis $\{U_{i}(0)\}_{i=1}^{n+m-1}$ is given by
\[
    A U_{j}(0) = \sum_{i=1}^{n+m-1} a_{ij} U_{i}(0).
\]

We now decompose equation \eqref{equ-Jacobi-equation} into its horizontal and vertical components:
\[
    \zeta_{j}^{h}{}'' + R_{c}^{h}(\gamma'^{h}, \zeta_{j}^{h})\gamma'^{h} = 0,
    \qquad
    \zeta_{j}^{v}{}'' + R_{c}^{v}(\gamma'^{v}, \zeta_{j}^{v})\gamma'^{v} = 0.
\]
Using the known solutions of Jacobi fields in $M_{c}^{n}$ and $\mathbb{R}^{m}$, we obtain
\begin{equation}\label{equ-What-bij(r)}
    \begin{cases}
        b_{ij}(r) = \delta_{ij} - a_{ij}r, & i \le m, \\[0.4em]
        b_{ij}(r) = \delta_{ij} C_{\tau}(r) - a_{ij} S_{\tau}(r), & i > m,
    \end{cases}
\end{equation}
where $S_{\tau}(r)$ and $C_{\tau}(r)$ are defined by
\[
    S_{\tau}(r) :=
    \begin{cases}
        \dfrac{1}{\sqrt{c\tau^{2}}}\sin\!\left(\sqrt{c\tau^{2}}\,r\right), & c > 0, \\[0.8em]
        \dfrac{1}{\sqrt{-c\tau^{2}}}\sinh\!\left(\sqrt{-c\tau^{2}}\,r\right), & c < 0,
    \end{cases}
    \quad
    C_{\tau}(r) :=
    \begin{cases}
        \cos\!\left(\sqrt{c\tau^{2}}\,r\right), & c > 0, \\[0.4em]
        \cosh\!\left(\sqrt{-c\tau^{2}}\,r\right), & c < 0.
    \end{cases}
\]
Moreover, these functions satisfy the first-order differential relations
\begin{equation}\label{equ-S'andC'}
    S_{\tau}'(r) = C_{\tau}(r),
    \qquad
    C_{\tau}'(r) = -c\tau^{2} S_{\tau}(r).
\end{equation}

In fact, the matrix $B(r) = (b_{ij}(r))$ given in equation \eqref{equ-What-bij(r)} can be written as the block matrix
\begin{align}\label{equ-B(r)-formula}
    B(r)
    =\left(
        \begin{array}{c|c}
            \delta_{ij}-a_{ij}r & -a_{ij}r \\
            \hline
            -a_{ij} S_{\tau}(r) & \delta_{ij} C_{\tau}(r) -a_{ij} S_{\tau}(r)
        \end{array}
    \right).
\end{align}

By Jacobi field theory, $B(r)$ is nonsingular for all $r \in (-\delta, \delta)$, and the shape operator of $\Sigma_{r}$ is given by  
\[
    A_{r} = -B'(r)B(r)^{-1}
    \quad \text{(cf. \cite[Theorem 10.2.1]{Berndt-Console-Olmos-2016})}.
\]
Hence, the mean curvature $H(r)$ is given by
\[
    H(r)
    = \operatorname{tr} A_{r}
    = -\operatorname{tr}\!\big(B'(r)B(r)^{-1}\big)
    = -\frac{d}{dr}(\det B(r)) / \det B(r).
\]
Defining $D(r) := \det B(r)$ and differentiating, we obtain
\[
    D'(r) + H(r)D(r) = 0,
\]
that is,
\[
    D'(r) = -H(r)D(r).
\]
By differentiating this equation repeatedly, for all $k \in \mathbb{N}$ we have
\begin{equation}\label{equ-0=Dk+1+phiD}
    0 = D^{(k+1)}(r) + \phi_{k}(r) D(r),
\end{equation}
where
\[
    \phi_{k}(r)
    = \phi_{k}\big(H(r), H'(r), \ldots, H^{(k)}(r)\big).
\]

\medskip

Recalling the structure of the matrix $B(r)$ in \eqref{equ-B(r)-formula}, we observe that the highest power of $r$ in the explicit expression for $D(r)$ is $m$. Hence, there exist coefficients $\alpha_{\ell,k}^{q}$ $(q=0,\ldots,m)$ such that
\begin{equation}\label{equ-Dk=Sumalpha_lk^q}
    D^{(k)}(r)
    =\sum_{\ell=0}^{n-1}\sum_{q=0}^{m}
    \alpha_{\ell,k}^{q} r^{q} S_{\tau}^{\ell}(r)C_{\tau}^{\,n-1-\ell}(r),
\end{equation}
where $D^{(k)}(r)$ denotes the $k$-th derivative of $D(r)$.

Substituting \eqref{equ-Dk=Sumalpha_lk^q} into \eqref{equ-0=Dk+1+phiD} and letting $k$ vary from $1$ to $(m+1)n-1$, we obtain
\begin{equation}\label{equ-0=alpha00k+1+phik}
    \alpha_{0,k+1}^{0}=-\phi_{k}(0).
\end{equation}

Using \eqref{equ-S'andC'} together with \eqref{equ-Dk=Sumalpha_lk^q}, we compute
\begin{align*}
    D^{(k+1)}(r)
    =&\sum_{\ell=0}^{n-1}\sum_{q=0}^{m}q\alpha_{\ell,k}^{q}r^{q-1}S_{\tau}^{\ell}(r)C_{\tau}^{n-1-\ell}(r)
    +\sum_{\ell=0}^{n-1}\sum_{q=0}^{m}\alpha_{\ell,k}^{q}r^{q}\ell S_{\tau}^{\ell-1}(r)C_{\tau}^{n-\ell}(r) \\
    &-\sum_{\ell=0}^{n-1}\sum_{q=0}^{m}\alpha_{\ell,k}^{q}r^{q}(n-1-\ell)c\tau^{2}S_{\tau}^{\ell+1}(r)C_{\tau}^{n-2-\ell}(r) \\
    =&\left(\sum_{q=0}^{m-1}\left(\left(q+1\right)\alpha_{0,k}^{q+1}+\alpha_{1,k}^{q}\right)r^{q}+\alpha_{1,k}^{m}r^{m}\right)C_{\tau}^{n-1}(r) \\
    &+\sum_{\ell=1}^{n-2}\left(\sum_{q=0}^{m-1}\left((q+1)\alpha_{\ell,k}^{q+1}+(\ell+1)\alpha_{\ell+1,k}^{q}
    -(n-\ell)c\tau^{2}\alpha_{\ell-1,k}^{q}\right)r^{q} \right. \\
    &\qquad\qquad \left. +\left((\ell+1)\alpha_{\ell+1,k}^{m}
    -(n-\ell)c\tau^{2}\alpha_{\ell-1,k}^{m}\right)r^{m}\right)
    S_{\tau}^{\ell}(r)C_{\tau}^{n-1-\ell}(r) \\
    &+\left(\sum_{q=0}^{m-1}((q+1)\alpha_{n-1,k}^{q+1}
    -c\tau^{2}\alpha_{n-2,k}^{q})r^{q}
    -c\tau^{2}\alpha_{n-2,k}^{m}r^{m}\right)
    S_{\tau}^{n-1}(r).
\end{align*}
Therefore, for $\ell=0,\ldots,n-1$ and $q=0,\ldots,m$, the coefficients satisfy
\begin{equation}\label{equ-alpha_k_to_alpha_k+1}
    \alpha_{\ell,k+1}^{q}
    =(q+1)\alpha_{\ell,k}^{q+1}
    +(\ell+1)\alpha_{\ell+1,k}^{q}
    -(n-\ell)c\tau^{2}\alpha_{\ell-1,k}^{q},
\end{equation}
where we set $\alpha_{\ell,k}^{m+1}=0$ for all $\ell=0,\ldots,n-1$ and
$\alpha_{-1,k}^{q}=\alpha_{n,k}^{q}=0$ for all $q=0,\ldots,m$.

From the recursive relation \eqref{equ-alpha_k_to_alpha_k+1} among the coefficients $\alpha_{\ell,k}^{q}$, we may write
\begin{equation}\label{equ-alpha00k+1=Sum_pqk+1l_alphapl0}
    \alpha_{0,k+1}^{0}
    =\sum_{\ell=0}^{n-1}\sum_{q=0}^{m}
     p_{k+1,\ell}^{q}\,\alpha_{\ell,0}^{q},
\end{equation}
where each coefficient $p_{k+1,\ell}^{q}$ depends only on the parameters
$q,k,\ell,n,m,c,$ and $\tau$.

Since $\alpha_{0,0}^{0}=D(0)=1$ and $\alpha_{0,k+1}^{0}$ coincides with $\phi_{k}(0)$ in equation \eqref{equ-0=alpha00k+1+phik}, we conclude that the vector
\begin{equation*}
    \xi_{0}
    =(\alpha_{1,0}^{0},\ldots,\alpha_{n-1,0}^{0},
      \alpha_{0,0}^{1},\ldots,\alpha_{n-1,0}^{m-1},
      \alpha_{0,0}^{m},\ldots,\alpha_{n-1,0}^{m})^{\mathrm{T}}
    \in\mathbb{R}^{(m+1)n-1}
\end{equation*}
satisfies a linear system of the form $M\xi=\nu$, according to \eqref{equ-0=alpha00k+1+phik}, where
\begin{equation*}
    \nu
    =(-\phi_{1}(0)-p_{2,0}^{0},\ldots,
      -\phi_{(m+1)n-1}(0)-p_{(m+1)n,0}^{0})^{\mathrm{T}}
    \in\mathbb{R}^{(m+1)n-1}.
\end{equation*}

\medskip

In the following, we shall see that the matrix $M$ exhibits fundamentally different properties depending on whether $n$ is odd or even.  
For odd $n$, we further denote by $M^{s}$ the $((m+1)n-1)\times((m+1)n-1)$ matrix obtained from $\widetilde{M}^{s}$ in Proposition \ref{prop-n-even-independent-odd-ss-dependent}-(ii) by removing its first column.  
For any $n$, let $M_{\iota}$ and $M_{\iota}^{s}$ denote the matrices obtained by replacing the $\iota$-th column of $M$ and $M^{s}$, respectively, with the vector $\nu$.

We establish key properties of $M$ and $M^{s}$ (for $s\ge (m+1)n$), in particular deriving a non-trivial algebraic expression in $\tau$.
Since the full proof is rather technical, it is deferred to the end of this section.
 \vspace{1mm}
 
\begin{Proposition}\label{Proposition-M-and-Ms}
    The matrices $M$ \((\)for even $n$\()\) and $M^{s}$ \((\)for odd $n$ and any $s\ge (m+1)n$\()\) satisfy the following properties:
    \begin{enumerate}[label=\textup{(\roman*)}]
        \item $\mathrm{rank}\,M=(m+1)n-2$ and $\mathrm{rank}\,M^{s}=(m+1)n-2$;
        \item There exists $\iota\in \{1,\ldots,(m+1)n-1\}$ such that
            \begin{equation*}
                \det M_{\iota}
                =(-1)^{\frac{\gamma_{0}}{2}}\beta_{0}c^{\frac{\gamma_{0}}{2}}\tau^{\gamma_{0}}
                -\sum_{i=1}^{(m+1)n-1}
                (-1)^{\frac{{\gamma}_{i}}{2}}\beta_{i}\phi_{i}(0)c^{\frac{\gamma_{i}}{{2}}}\tau^{\gamma_{i}},
            \end{equation*}
            where $\beta_{0}\neq 0$, and $\beta_{1},\ldots,\beta_{(m+1)n-1}$ as well as $\gamma_{0}>\cdots>\gamma_{(m+1)n-1}>0$ are integers;
        \item 
            \begin{equation*}
                \det M_{n}^{s}
                =-(-1)^{\frac{\gamma_{s}}{2}}\beta_{s}\phi_{s}(0)c^{\frac{\gamma_{s}}{2}}\tau^{\gamma_{s}}
                -\sum_{i=1}^{(m+1)n-2}
                (-1)^{\frac{\gamma_{i}}{2}}\beta_{i}\phi_{i}(0)c^{\frac{\gamma_{i}}{2}}\tau^{\gamma_{i}},
            \end{equation*}
            where $\beta_{s}\neq 0$,
            $\beta_{1},\ldots,\beta_{(m+1)n-2}$ as well as $\gamma_{1}>\cdots >\gamma_{(m+1)n-2}>\gamma_{s}>0$ are all integers.
    \end{enumerate}
\end{Proposition}
 \vspace{1mm}
 
\noindent
\textit{Proof of Theorem \ref{thm-main}.$\quad$}
We will primarily apply the non-trivial algebraic expression in $\tau$ implied by $M$ and $M^{s}$ as established in Proposition \ref{Proposition-M-and-Ms}, and prove Theorem \ref{thm-main} proceeding case by case.

\medskip

\noindent\textbf{Case 1:} $n \ge 2$ and $n$ is even.  
By Proposition \ref{Proposition-M-and-Ms}-(i), $\det M = 0$.
Since $\xi_{0}$ satisfies $M\xi = \nu$, it follows that $\det M_{i} = 0$ for all $i = 1, \ldots, (m+1)n-1$.  
Moreover, by Proposition \ref{Proposition-M-and-Ms}-(ii), there exists an index $\iota \in \{1,\ldots, (m+1)n-1\}$ such that
\begin{align*}
    \det M_{\iota}
    =(-1)^{\frac{\gamma_{0}}{2}}\beta_{0}c^{\frac{\gamma_{0}}{2}}\tau^{\gamma_{0}}
    -\sum_{i=1}^{(m+1)n-1} (-1)^{\frac{\gamma_{i}}{2}}\beta_{i}\phi_{i}(0)c^{\frac{\gamma_{i}}{2}}\tau^{\gamma_{i}}
    = 0.
\end{align*}
This yields a nontrivial algebraic equation in $\tau$, and hence $\tau$ is constant.

\medskip

\noindent\textbf{Case 2:} $n \ge 2$ and $n$ is odd.  
If there exists $s_{0} \ge (m+1)n$ with $\phi_{s_{0}}(0) \neq 0$, then, analogous to Case 1, Proposition \ref{Proposition-M-and-Ms} (i) and (iii) imply that
\[
    \det M^{s_{0}}_{n} = 0
\]
is a nontrivial algebraic equation in $\tau$, so $\tau$ is constant.  

Otherwise, if $\phi_{s}(0) = 0$ for all $s \ge (m+1)n$, then equation \eqref{equ-0=Dk+1+phiD} shows that $D(r)$ is polynomial near $r=0$.  
However, from equation \eqref{equ-B(r)-formula}, $D(r) = \det B(r)$ cannot be polynomial near $r=0$, so this case is excluded.

\hfill$\Box$


\subsection{Proof of Proposition \ref{Proposition-M-and-Ms}}

We first derive the recurrence relation in $k$ for the coefficients $p_{k+1,\ell}^{q}$ in equation \eqref{equ-alpha00k+1=Sum_pqk+1l_alphapl0}.  
For $n \ge 4$, combining equations \eqref{equ-alpha_k_to_alpha_k+1} and \eqref{equ-alpha00k+1=Sum_pqk+1l_alphapl0}, we obtain the following explicit computation:
\begin{align*}
    \alpha_{0,k+1}^{0}
    =&\sum_{\ell=0}^{n-1}\sum_{q=0}^{m}p_{k,\ell}^{q}\alpha_{\ell,1}^{q} \\
    =&\sum_{q=0}^{m}p_{k,0}^{q}\alpha_{0,1}^{q}
    +\sum_{\ell=1}^{n-2}\sum_{q=0}^{m}p_{k,\ell}^{q}\alpha_{\ell,1}^{q}
    +\sum_{q=0}^{m}p_{k,n-1}^{q}\alpha_{n-1,1}^{q} \\
    =&\sum_{q=1}^{m}p_{k,0}^{q-1}q\alpha_{0,0}^{q}
    +\sum_{q=0}^{m-1}p_{k,0}^{q}\alpha_{1,0}^{q}
    +p_{k,0}^{m}\alpha_{1,0}^{m} \\
    &+\sum_{\ell=1}^{n-2}\sum_{q=1}^{m}p_{k,\ell}^{q-1}q\alpha_{\ell,0}^{q}
    +\sum_{\ell=2}^{n-1}\sum_{q=0}^{m-1}p_{k,\ell-1}^{q}\ell\alpha_{\ell,0}^{q}
    -\sum_{\ell=0}^{n-3}\sum_{q=0}^{m-1}p_{k,\ell+1}^{q}(n-1-\ell)c\tau^{2}\alpha_{\ell,0}^{q} \\
    &+\sum_{\ell=2}^{n-1}p_{k,\ell-1}^{m}\ell\alpha_{\ell,0}^{m}
    -\sum_{\ell=0}^{n-3}p_{k,\ell+1}^{m}(n-1-\ell)c\tau^{2}\alpha_{\ell,0}^{m} \\
    &+\sum_{q=1}^{m}p_{k,n-1}^{q-1}q\alpha_{n-1,0}^{q}
    -\sum_{q=0}^{m-1}p_{k,n-1}^{q}c\tau^{2}\alpha_{n-2,0}^{q}
    -p_{k,n-1}^{m}c\tau^{2}\alpha_{n-2,0}^{m} \\
    =&\sum_{\ell=0}^{n-1}\sum_{q=0}^{m}
    \left(
    p_{k,\ell}^{q-1}q
    +p_{k,\ell-1}^{q}\ell
    -p_{k,\ell+1}^{q}(n-1-\ell)c\tau^{2}
    \right) \alpha_{\ell,0}^{q}.
\end{align*}
For convenience, setting $p_{k,\ell}^{q}=0$ if $q=-1$ or $\ell =-1, n$, we have
\begin{equation}\label{equ-Relpqk+1landpqkl}
    p_{k+1,\ell}^{q}
    =q\, p_{k,\ell}^{q-1}
    +\ell\, p_{k,\ell-1}^{q}
    -(n-1-\ell)c\tau^{2}p_{k,\ell+1}^{q},
\end{equation}
for any $q=0,\ldots, m$ and $\ell=0,\ldots, n-1$.
 \vspace{1mm}
 
\begin{Remark}
For $n=2$ and $n=3$,
the computations are entirely analogous and remain relatively straightforward.
Although we omit the explicit calculations for brevity, these low-dimensional cases reproduce the formula in equation \eqref{equ-Relpqk+1landpqkl} exactly.  
This confirms that equation \eqref{equ-Relpqk+1landpqkl} holds for all integers $n \ge 2$.
\end{Remark}
 \vspace{1mm}
 
More specifically, we have $p_{0,0}^{0}=1$ and $p_{0,\ell}^{q}=0$ for all other $(\ell,q)$, due to $\alpha_{0,0}^{0}=1$ and equation \eqref{equ-alpha00k+1=Sum_pqk+1l_alphapl0}.
Taking $k=0$ in \eqref{equ-Relpqk+1landpqkl}, we obtain
\begin{equation*}
    p_{1,0}^{1} = p_{1,1}^{0} = 1,
    \quad
    p_{1,\ell}^{q} = 0
    \quad 
    \text{for } (\ell,q) \neq (0,1), (1,0).
\end{equation*}
Similarly, when $k=1$, we find that
\begin{equation}\label{equ-equ-Relpqk+1landpqkl-case-p2}
    p_{2,\ell}^{q} =
    \begin{cases}
        -(n-1)c\tau^{2}, & \text{if~}(\ell,q)=(0,0), \\
        2, & \text{if~}(\ell,q)=(0,2),(1,1),(2,0), \\
        0, & \text{otherwise}.
    \end{cases}
\end{equation}

By applying the recurrence relation \eqref{equ-Relpqk+1landpqkl} and mathematical induction, we obtain the following basic characterization of $p_{k,\ell}^{q}$:
 \vspace{1mm}
 
\begin{Proposition}\label{prop-how-pklq}
	When $n\ge 2$, for any $k\ge 2$, $q=0,\ldots,m$, $\ell=0,\ldots,n-1$, the identity $p_{k,\ell}^{q}=\sigma_{k,\ell}^{q}(n)c^{s}\tau^{2s}$ holds with $s=\frac{1}{2}(k-q-\ell)$.
	Furthermore, the following assertions also hold:
	\begin{enumerate}[label=\textup{(\roman*)}]
		\item $\sigma_{k,\ell}^{q}(n)=0$, for $s\not\in \mathbb{Z}$ or $s<0$. \label{prop-how-pklq-item-i}
		\item $\sigma_{k,\ell}^{q}(n)=k!$, for $s=0$. \label{prop-how-pklq-item-ii}
		\item $\sigma_{k,\ell}^{q}(n)$ is a polynomial of degree $\deg\sigma_{k,\ell}^{q}(n)\ge s$ with $(-1)^s$ as its leading coefficient sign for $s\in \mathbb{Z}_{+}$. \label{prop-how-pklq-item-iii}
	\end{enumerate}
\end{Proposition}

\begin{Proof}
	(i)
	If $s \notin \mathbb{Z}$, then $k - q - \ell$ is odd.
	As discussed above, $p_{0,0}^{0}=1$ and $p_{0,\ell}^{q}=0$ for all $\ell=1,\ldots,n-1$ and $q=0,\ldots,m$.  
	Since the parity of $k - q - \ell$ is preserved in equation~\eqref{equ-Relpqk+1landpqkl}, it follows that $p_{k,\ell}^{q}=0$ for all $s \notin \mathbb{Z}$.
	
	If $s < 0$, then $q + \ell > k$. By induction, we will show that $p_{k,\ell}^{q} = 0$ also holds in this case.
	
	For $k = 2$, by equation \eqref{equ-equ-Relpqk+1landpqkl-case-p2}, we have  
	\[
	p_{2,\ell}^{q} = 0, \quad \text{for all $q, \ell$ with $q + \ell > 2$}.
	\]
	
	Now assume that for some $k \ge 2$, $p_{k,\ell}^{q} = 0$ for all $(q, \ell)$ satisfying $q + \ell > k$.  
	Then for any $q, \ell$ with $q + \ell > k + 1$, equation \eqref{equ-Relpqk+1landpqkl} gives
	\[
	p_{k+1,\ell}^{q}
	= q\,p_{k,\ell}^{q-1}
	+ \ell\,p_{k,\ell-1}^{q}
	- (n - 1 - \ell)c\tau^{2}p_{k,\ell+1}^{q}.
	\]
	By the induction hypothesis, $p_{k,\ell}^{q-1} = p_{k,\ell-1}^{q} = p_{k,\ell+1}^{q} = 0$ since $(q - 1) + \ell > k$, $q + (\ell - 1)>k$, and $q + (\ell + 1)> k$.  
	Therefore, $p_{k+1,\ell}^{q} = 0$, as required.

	(ii)
	It suffices to show that $p_{k,k-q}^{q} = k!$ for all such $k,q$.
	When $k=2$, equation \eqref{equ-equ-Relpqk+1landpqkl-case-p2} gives $p_{2,2}^{0}=2$.
	
	Suppose that for some integer $k\geq 2$, the identity $p_{k,k-q}^{q} = k!$ holds for all $q$. Then, for any $q$ satisfying $k + 1 \geq 2 + q$,
	equation \eqref{equ-Relpqk+1landpqkl} implies 
	\begin{align*}
		p_{k+1,k+1-q}^{q} 
		&= q\, p_{k,k+1-q}^{q-1}
		+ (k+1-q)\, p_{k,k-q}^{q}
		- (n-1-(k+1-q)) c \tau^{2} p_{k,k+2-q}^{q} \\
		&= (k+1)!,
	\end{align*}
	where the last equality uses $p_{k,k+1-q}^{q-1} = p_{k,k-q}^{q} = k!$ from the induction hypothesis and $p_{k,k+2-q}^{q} = 0$ from (i).
	Hence, the result follows by induction.
	
	(iii)
	Equivalently, it suffices to show that for any \(k \ge 2\),
	\[
	\deg \sigma_{k,\ell}^{q}(n) \ge s, \qquad \text{where } \ell = k - q - 2s,
	\]
	for all integers \(s\) and \(q\) satisfying \(0 < s \le \tfrac{1}{2}(k - q)\).
	
	To prove this, we start with \(k=2\) and proceed by induction.  
	From equation \eqref{equ-equ-Relpqk+1landpqkl-case-p2}, the only term satisfying the condition is \(p_{2,0}^{0} = -(n-1)c\tau^{2}\).  
	Hence \(\deg \sigma_{2,0}^{0}(n) = s = 1\), verifying the claim for \(k=2\).
	
	Assume now that for some \(k \ge 2\), one has \(\deg \sigma_{k,\,k-q-2s}^{q}(n) \ge s\) for all integers \(q, s\) satisfying \(0 < s \le \tfrac{1}{2}(k - q)\), and that the leading coefficient has sign \( (-1)^{s} \).
	For \(k+1\), take any such \(q, s\) with \(0 < s \le \tfrac{1}{2}(k+1 - q)\) and set \(\ell = k + 1 - q - 2s\).
	By equation \eqref{equ-Relpqk+1landpqkl}, we obtain
	\begin{align*}
		p_{k+1,\ell}^{q} 
		&= q\,p_{k,\ell}^{q-1} + \ell\,p_{k,\ell-1}^{q}
		- (n-1-\ell)c\tau^{2}p_{k,\ell+1}^{q} \\
		&= q\, \sigma_{k,\ell}^{q-1}(n) c^{s} \tau^{2s} 
		+ \ell\, \sigma_{k,\ell-1}^{q}(n) c^{s} \tau^{2s} 
		- (n-1 - \ell) \sigma_{k,\ell+1}^{q}(n) c^{s} \tau^{2s} \\
		&= \sigma_{k+1,\ell}^{q}(n)c^{s}\tau^{2s},
	\end{align*}
	where
	\[
	\sigma_{k+1,\ell}^{q}(n)
	= q\,\sigma_{k,\ell}^{q-1}(n)
	+ \ell\,\sigma_{k,\ell-1}^{q}(n)
	- (n-1-\ell)\sigma_{k,\ell+1}^{q}(n).
	\]
	By the induction hypothesis, we have 
	\(\deg \sigma_{k,\ell}^{q-1}(n) = \deg \sigma_{k,\ell-1}^{q} = s\) and \(\deg \sigma_{k,\ell+1}^{q} = s-1\), 
	with respective leading coefficient signs \( (-1)^{s}, (-1)^{s} \), and \( (-1)^{s-1} \). 
	It then follows that \(\deg \sigma_{k+1,\ell}^{q}(n) \ge s\) and its leading coefficient has sign \( (-1)^{s} \), 
	which completes the proof.
	
\end{Proof}

To study the rank properties of $M$ and $M^{s}$ ($s \ge (m+1)n$), we adopt a row-wise perspective.  
Define the matrix $ \widetilde{M} = \left[ -\nu_{\tau},\, M \right] $ where $\nu_{\tau} = \nu - \nu_{\phi}$ and
\begin{equation*}
    \nu_{\phi}=(-\phi_{1}(0),\ldots,-\phi_{(m+1)n-1}(0))^{\mathrm{T}}.
\end{equation*} 
Let $M_\iota^\tau$ (resp.,\ $M_\iota^\phi$) denote the matrix obtained from $M$ by replacing its $\iota$-th column with $\nu_\tau$ (resp.,\ $\nu_\phi$).
Under this setting,
each row of the matrix $\widetilde{M}$ is a row vector of the form
\[
\widetilde{L}_{k-1} = \big(p_{k,0}^{0}, \ldots, p_{k,n-1}^{0},\, p_{k,0}^{1}, \ldots, p_{k,n-1}^{1},\, \ldots,\, p_{k,0}^{m}, \ldots, p_{k,n-1}^{m} \big), \quad k \ge 2,
\]
where each segment $p_{k,0}^{q},\dots,p_{k,n-1}^{q}$ corresponds to $q=0,\dots,m$.  

Define $e_{1} \in \mathbb{R}^n$ by
\[
e_{1} = (p_{0,0}^{0}, p_{0,1}^{0}, \dots, p_{0,n-1}^{0}) = (1,0,\dots,0).
\]  

Next, we define an $(m+1)n \times (m+1)n$ matrix $Q$ as follows:
\begin{equation}\label{equ-What-Q}
Q =
\begin{pmatrix}
K & I & & & & & \\
& K & 2I & & & & \\
& & K & 3I & & & \\
& & & \ddots & \ddots & & \\
& & & & K & (m-1)I & \\
& & & & & K & mI \\
& & & & & & K
\end{pmatrix},
\end{equation}
where $I$ is the $n \times n$ identity matrix and $K$ is the $n \times n$ $\tau$-Kac matrix
\[
K =
\begin{pmatrix}
0 & 1 & 0 & \cdots & 0 & 0 & 0 \\
-(n-1)c\tau^2 & 0 & 2 & \cdots & 0 & 0 & 0 \\
0 & -(n-2)c\tau^2 & 0 & \cdots & 0 & 0 & 0 \\
\vdots & \vdots & \vdots & \ddots & \vdots & \vdots & \vdots \\
0 & 0 & 0 & \cdots & 0 & n-2 & 0 \\
0 & 0 & 0 & \cdots & -2c\tau^2 & 0 & n-1 \\
0 & 0 & 0 & \cdots & 0 & -c\tau^2 & 0
\end{pmatrix}.
\]

Then, for $k \ge 2$, from the recurrence equation \eqref{equ-Relpqk+1landpqkl}, the rows of $\widetilde{M}$ satisfy
\[
\widetilde{L}_k = \widetilde{L}_{k-1} Q = (e_1,0,\dots,0) Q^{k+1}.
\]

Regarding the properties of the $ \tau $-Kac matrix,
we recall the following lemma established in \cite{de-Lima-Pipoli-2024,Edelman-Kostlan-1994}.
 \vspace{1mm}
 
\begin{Lemma}\textup{(\cite{de-Lima-Pipoli-2024,Edelman-Kostlan-1994})}\label{lem-de-Lima-Lemma16}
    The $\tau$-Kac matrix $K$ of order $n$ has the following properties:
    \begin{enumerate}[label=\textup{(\roman*)}]
        \item It has $n$ simple eigenvalues $\lambda_{0},\ldots,\lambda_{n-1}$,
        which are
        \begin{equation*}
            \lambda_{\ell}
            =(n-1-2\ell)\sqrt{-c}\tau,
            \quad
            \ell\in\{0,\ldots,n-1\}.
        \end{equation*}
        In particular $\lambda_{\ell}$ is real if $c<0$,
        and purely imaginarg if $c>0$;
        \item Its rank is $n$,
        if $n$ is even,
        and $n-1$ if $n$ is odd.
        In particular,
        $K$ is nonsignular if and only if $n$ is even;
        \item The coordinates of $e_{1}\in\mathbb{R}^{n}$ with respect to the basis of its eigenvectors are all different from zero.
    \end{enumerate}
\end{Lemma}
 \vspace{1mm}
 
We now relate $Q$ to the $\tau$-Kac matrix $K$, in particular its eigenvectors.
Direct computations give
\[
\det Q = (\det K)^{m+1},
\]
and thus by Lemma \ref{lem-de-Lima-Lemma16}-(ii), $Q$ is nonsingular if and only if $n$ is even.

Let $\{x_0, \ldots, x_{n-1}\} \subset \mathbb{R}^n$ be the eigenvectors of $K$.  
For $\ell = 0, \ldots, n-1$, we define the following vector in $ \mathbb{R}^{(m + 1)n} $:
\[
x_{0,\ell} = (x_\ell, 0, \ldots, 0),\quad
x_{1,\ell} = (0, x_\ell, 0, \ldots, 0),\quad \dots,\quad
x_{m,\ell} = (0, \ldots, 0, x_\ell).
\]
Then, we have
\[
x_{k,\ell} Q = \lambda_\ell x_{k,\ell} + x_{k+1,\ell}, \quad k=0,\ldots,m-1, 
\qquad
x_{m,\ell} Q = \lambda_\ell x_{m,\ell}.
\]
More generally, for any integer $k \ge 0$ and $0 \le i \le m$, by induction,
\begin{equation}\label{equ-xiQk=xi}
x_{i,\ell} Q^k = \sum_{t=0}^{\min\{k, m-i\}} \binom{k}{t} \lambda_\ell^{\,k-t} x_{i+t,\ell}.
\end{equation}

 \vspace{1mm}
\begin{Proposition}\label{prop-n-even-independent-odd-ss-dependent}
Let $\widetilde{e}_{1}=(e_{1},0,\ldots,0)\in\mathbb{R}^{(m+1)n}$ with $n\ge 2$. Then, the following assertions hold.
\begin{enumerate}[label=\textup{(\roman*)}]
    \item If $n$ is even, for any positive integer $s$, the set
    \[
        \{\widetilde{e}_{1} Q^{i} \mid i=s, \ldots, s+(m+1)n-1\}
    \]
    is linearly independent.

    \item If $n$ is odd, for any integer $s \ge (m+1)n$, define the ordered sets
    \[
        \Lambda = \{\widetilde{e}_{1} Q^{i} \mid i = 2, \ldots, (m+1)n-1\}, 
        \qquad
        \Lambda_{s} = \Lambda \cup \{\widetilde{e}_{1} Q^{s}\}.
    \]
    Let $\widetilde{M}^{s}$ be the matrix with rows given by the vectors in $\Lambda_{s}$, and denote its columns by $C_{1},\ldots,C_{(m+1)n}$. Then the following hold:
    \begin{enumerate}[label=\textup{(\alph*)}]
        \item $\Lambda$ is linearly independent, whereas $\Lambda_{s}$ is linearly dependent.
        \item For $q = 0, 1$, the column $C_{qn+1}$ lies in the span of the columns $\{C_{qn+2i+1} \mid i = 1, \ldots, (n-1)/2\}$.
    \end{enumerate}
\end{enumerate}
\end{Proposition}
 \vspace{1mm}
 
\begin{Proof}
(i) When $n$ is even, the previous calculations show that $Q$ is invertible, so it suffices to consider $s=0$.

Consider the vector equation
\begin{equation}\label{equ-e1Q_independent-i-1}        
    \sum_{k=0}^{(m+1)n-1} \mu_{k} \widetilde{e}_{1} Q^{k} = 0    
\end{equation}
in the variables $\mu_{0}, \ldots, \mu_{(m+1)n-1}$. 

Without loss of generality, by Lemma \ref{lem-de-Lima-Lemma16}-(iii), we may write
\begin{equation*}        
    \widetilde{e}_{1} = \sum_{\ell=0}^{n-1} a_{\ell} x_{0,\ell},
\end{equation*}
where $a_{\ell} \neq 0$ for all $\ell = 0, \ldots, n-1$.

From equation \eqref{equ-xiQk=xi}, we obtain
\begin{equation*}
    \sum_{\ell=0}^{n-1}\sum_{k=0}^{(m+1)n-1}\sum_{t=0}^{\min\{k,m\}}\binom{k}{t}\lambda_{\ell}^{k-t}\mu_{k}a_{\ell}x_{t,\ell}=0.
\end{equation*}
Hence, the system \eqref{equ-e1Q_independent-i-1} is equivalent to the linear system
\begin{equation}\label{equ-e1Q_independent-i-2}        
    \sum_{k=0}^{(m+1)n-1} \binom{k}{t} \lambda_{\ell}^{\,k-t} \mu_{k} = 0, \quad t = 0,\ldots,m,\;\ell=0,\ldots,n-1.
\end{equation}

The coefficient matrix $\Xi$ of \eqref{equ-e1Q_independent-i-2} is a generalized Vandermonde matrix with
\begin{equation*}        
    \det \Xi = \prod_{i<j} (\lambda_j - \lambda_i)^{(m+1)^2}.
\end{equation*}
By Lemma \ref{lem-de-Lima-Lemma16}, the eigenvalues $\lambda_\ell$ are distinct, so $\det \Xi \neq 0$ and $\Xi$ is invertible. Hence, $\mu_{k} = 0$ for all $k = 0, \ldots, (m+1)n-1$, completing the proof of (i).

(ii)-(a) Similar to (i), for any $s \ge (m+1)n$, we consider the following vector equation in the variables $\mu_{2},\ldots,\mu_{(m+1)n-1},\mu_{s}$:
\[
\sum_{k=2}^{(m+1)n-1} \mu_k \widetilde{e}_1 Q^k + \mu_s \widetilde{e}_1 Q^s = 0,
\]
which is equivalent to the linear system
\begin{equation}\label{equ-equ-e1Q_independent-ii-a-2}
\sum_{k=2}^{(m+1)n-1} \binom{k}{t} \lambda_\ell^{\,k-t} \mu_k
+ \binom{s}{t} \lambda_\ell^{\,s-t} \mu_s = 0, \quad t = 0, \ldots, m,\;\ell = 0, \ldots, n-1.
\end{equation}

Since $n$ is odd, Lemma \ref{lem-de-Lima-Lemma16}-(i) gives $\lambda_{(n-1)/2}=0$. Thus, \eqref{equ-equ-e1Q_independent-ii-a-2} is a linear system of $(m+1)n-2$ equations in $(m+1)n-1$ unknowns. Its coefficient matrix $\Xi$ has the block form
\[
\Xi =
\begin{pmatrix}
\Xi_0 \\ \Xi_1 \\ \vdots \\ \Xi_{\frac{n-3}{2}} \\ \Xi_{\frac{n-1}{2}} \\ \Xi_{\frac{n+1}{2}} \\ \vdots \\ \Xi_{n-1}
\end{pmatrix},
\]
where the block $\Xi_\ell$ is the $(m+1) \times ((m+1)n-1)$ matrix
\[
\Xi_{\ell}
=\begin{pmatrix}
    \lambda_{\ell}^{2} & \lambda_{\ell}^{3} & \lambda_{\ell}^{4} & \cdots & \lambda_{\ell}^{(m+1)n-1} & \lambda_{\ell}^{s} \\
    2\lambda_{\ell} & 3\lambda_{\ell}^{2} & 4\lambda_{\ell}^{3} & \cdots & ((m+1)n-1)\lambda_{\ell}^{(m+1)n-2} & s\lambda_{\ell}^{s-1} \\
    1 & 3\lambda_{\ell} & 6\lambda_{\ell}^{2} & \cdots & \binom{(m+1)n-1}{2} \lambda_{\ell}^{(m+1)n-3} & \binom{s}{2}\lambda_{\ell}^{s-2} \\
    \vdots & \vdots & \vdots & \ddots & \vdots & \vdots \\
    0 & 0 & 0 & \cdots & \binom{(m+1)n-1}{m}\lambda_{\ell}^{(m+1)(n-1)} & \binom{s}{m} \lambda_{\ell}^{s-m}
\end{pmatrix}.
\]
In particular,
\begin{equation}\label{equ-block-Xi-(n-1)/2}
\Xi_{\frac{n-1}{2}} =
\begin{pmatrix}
\mathrm{O} & \mathrm{O} \\
I_{m-1} & \mathrm{O}
\end{pmatrix}.
\end{equation}
Hence, $\mathrm{rank}\,\Xi \le (m+1)n-2$, implying that the set $\Lambda_s$ is linearly dependent.
\vspace{1mm}

To prove that $\Lambda$ is independent, remove the last column of $\Xi$ (corresponding to $\mu_s$) to form $\widetilde{\Xi}$. Expanding the block $\widetilde{\Xi}_{\frac{n-1}{2}}$ in equation \eqref{equ-block-Xi-(n-1)/2} and applying generalized Vandermonde determinant properties, we obtain
\[
\det \widetilde{\Xi} = \prod_{\substack{i<j\\ i,j \neq (n-1)/2}} (\lambda_j - \lambda_i)^{(m+1)^2} \cdot \prod_{i \neq (n-1)/2} \lambda_i^{(m+1)^2} \neq 0.
\]
Hence, $\widetilde{\Xi}$ is nonsingular, and $\Lambda$ is linearly independent.

(ii)-(b) By equation \eqref{equ-What-Q} and induction, we have
\begin{equation}\label{equ-What-Q^j-version1}
    Q^{j}
    =\begin{pmatrix}
        K^{j} & \dbinom{j}{1}1^{\overline{1}} K^{j-1} & \dbinom{j}{2}1^{\overline{2}} K^{j-2} 
        & \cdots 
        & \dbinom{j}{m}1^{\overline{m}} K^{j-m} \\
        0 & K^j & \dbinom{j}{1}2^{\overline{1}} K^{j-1} 
        & \cdots  
        & \dbinom{j}{m-1}2^{\overline{m-1}} K^{j-m+1} \\
        \vdots & \vdots & \vdots 
        & \ddots 
        & \vdots \\
        0 & 0 & 0 
        & \cdots 
        & \dbinom{j}{1}m^{\overline{1}} K^{j-1} \\
        0 & 0 & 0 
        & \cdots 
        & K^{j}
    \end{pmatrix},
\end{equation}
or equivalently,
\begin{equation*}
    Q^j[p,q] =
    \begin{cases}
    \dbinom{j}{d} p^{\overline{d}} K^{j-d}, & q = p + d, \; 1 \le d \le m, \\
    K^j, & q = p, \\
    0, & \text{otherwise},
    \end{cases}
\end{equation*}
where $p^{\overline{d}} = p(p+1)\cdots(p+d-1)$ is the rising factorial and $Q^{j}[p,q]$ is the element in the $p$-th row and $q$-th column of the matrix $Q^{j}$.

Similar to the above argument in (ii)-(a),
the last row will be immaterial.
Without loss of generality, assume $s = (m+1)n$. By equation \eqref{equ-What-Q^j-version1},
we note that $C_{qn+1},C_{qn+2},\ldots,C_{(q+1)n}$ are the columns of the matrix whose rows are
\[
e_1 \binom{2}{q} K^{2-q},\, e_1 \binom{3}{q} K^{3-q},\, \dots,\, e_1 \binom{(m+1)n}{q} K^{(m+1)n-q},
\]
for $q=0,1$.

We claim that the set $\{C_{qn+2i+1} \mid i = 0, \dots, (n-1)/2\}$ spans a space of dimension $(n-1)/2$. Consider the vector equation in $\bar{n} := \big\lfloor \frac{(m+1)n-q}{2} \big\rfloor$ variables $\mu_1,\dots,\mu_{\bar{n}}$:
\[
\sum_{j=1}^{\bar{n}} \mu_j \binom{2j+q}{q} e_1 K^{2j} = 0,
\]
which is equivalent to the linear system
\begin{equation}\label{equ-equ-e1Q_independent-ii-b-1}
    \sum_{j=1}^{\bar{n}} \binom{2j+q}{q} \lambda_\ell^{2j} \mu_j = 0, \quad \ell = 0, \dots, n-1.
\end{equation}

The coefficient matrix of the linear system \eqref{equ-equ-e1Q_independent-ii-b-1} is
\[
\Xi
=\begin{pmatrix}
    \dbinom{2+q}{q}\lambda_{0}^{2} & \dbinom{4+q}{q}\lambda_{0}^{4} & \cdots & \dbinom{2\bar{n}+q}{q}\lambda_{0}^{2\bar{n}} \\
    \dbinom{2+q}{q}\lambda_{1}^{2} & \dbinom{4+q}{q}\lambda_{1}^{4} & \cdots & \dbinom{2\bar{n}+q}{q}\lambda_{1}^{2\bar{n}} \\
    \vdots & \vdots & \ddots & \vdots \\
    \dbinom{2+q}{q}\lambda_{n-1}^{2} & \dbinom{4+q}{q}\lambda_{n-1}^{4} & \cdots & \dbinom{2\bar{n}+q}{q}\lambda_{n-1}^{2\bar{n}}
\end{pmatrix}.
\]

By Lemma \ref{lem-de-Lima-Lemma16}-(i), we know that $\lambda_{(n-1)/2} = 0$ and the nonzero eigenvalues occur in pairs $\pm \lambda$, and hence $\mathrm{rank}\,\Xi \le (n-1)/2$. Taking the first $(n-1)/2$ columns and rows corresponding to distinct eigenvalues gives a submatrix $\widetilde{\Xi}$ with
\[
\det \widetilde{\Xi} = \prod_{i<j} (\lambda_j^2 - \lambda_i^2) \cdot \prod_{\ell=0}^{(n-3)/2} \lambda_\ell^2 \neq 0.
\]
Therefore,
we obtain that $\mathrm{rank}\,\Xi=\frac{n-1}{2}$ and the claim is proved.

Finally, by Proposition \ref{prop-how-pklq}-(i) and (ii), the submatrix consisting of the first $\frac{n-1}{2}$ nonzero rows of the matrix formed by the column vectors $\{C_{qn+2i+1} \mid i = 1, \ldots, \frac{n-1}{2}\}$ is a lower triangular matrix, whose diagonal entries are $(2+q)!, (4+q)!, \ldots, (q+n-1)!$. Therefore, the vectors $\{C_{qn+2i+1} \mid i = 1, \ldots, \frac{n-1}{2}\}$ are linearly independent, and $C_{qn+1}$ lies in their span. This completes the proof.

\end{Proof}

\noindent
\textit{Proof of Proposition \ref{Proposition-M-and-Ms}.$\quad$}
(i) When $n$ is even, by Proposition \ref{prop-how-pklq}-(i), we observe that all the odd rows of the matrix $M$ form a matrix with $\frac{(m+1)n}{2}$ rows but only $\frac{(m+1)n}{2} - 1$ nonzero columns. Hence, the rank of the odd rows is at most $\frac{(m+1)n}{2}-1$, implying
\[
\mathrm{rank}\, M \le (m+1)n - 2.
\]

On the other hand, applying Proposition \ref{prop-n-even-independent-odd-ss-dependent}-(i) with $s=2$, the augmented matrix $\widetilde{M}$ has rank $(m+1)n - 1$, which shows
\[
\mathrm{rank}\, M \ge (m+1)n - 2.
\]
Combining the bounds, we conclude that $\mathrm{rank}\, M = (m+1)n - 2$.

When $n$ is odd, consider $s=(m+1)n$ in Proposition \ref{prop-n-even-independent-odd-ss-dependent}-(ii). The augmented matrix $\widetilde{M}^s$ then has rank $(m+1)n-2$. Moreover, by Proposition \ref{prop-n-even-independent-odd-ss-dependent}-(ii)-(b) with $q=0$, the first column $C_1$ lies in the span of $C_3, C_5, \dots, C_n$. Therefore, 
\[
\mathrm{rank}\, M^s = \mathrm{rank}\, \widetilde{M}^s = (m+1)n - 2.
\]

(ii)
By Proposition \ref{prop-n-even-independent-odd-ss-dependent}-(i), the augmented matrix $\widetilde{M}$ has rank $(m+1)n-1$. Hence, there exists an index $\iota \in \{1, \ldots, (m+1)n-1\}$ such that 
\[
\det M^{\tau}_{\iota} \neq 0.
\] 
By Proposition \ref{prop-how-pklq}-(iii), there exist a nonzero integer $\beta_0$ and a positive integer $\gamma_0 > 0$ such that
\[
\det M^{\tau}_{\iota} = (-1)^{\gamma_0/2} \beta_0 \, c^{\gamma_0/2} \tau^{\gamma_0}.
\]

Next, performing a Laplace expansion along the $\iota$-th column of $\det M^{\phi}_{\iota}$, we obtain the form described in (ii), with constant coefficients $\{\beta_k\}_{k=1}^{(m+1)n-1}$ and exponents $\{\gamma_k\}_{k=1}^{(m+1)n-1}$.  
Meanwhile, Proposition \ref{prop-how-pklq} implies that all $\gamma_{i}$ are even, each $\beta_{i}$ is an integer, and that the sequence of $\gamma_{i}$ is strictly increasing in $i$.

Finally, the positivity of $\gamma_{(m+1)n-1}$ follows from the estimate:
\begin{align*}
\gamma_{(m+1)n-1} 
&\ge \sum_{i=1}^{(m+1)n-2} (i-1) + \sum_{q=0}^{m} \sum_{j=1}^{n-1} \big(q(n-1) - j + 2\big) \\
&= \frac{1}{2}\Big( m^2 (2n^2 - 2n + 1) + m(2n^2 - 2n - 3) \Big) + 1 > 0.
\end{align*}

(iii) 
Define the matrices $M_{\iota}^{s}$, $M_{\iota}^{s,\tau}$, and $M_{\iota}^{s,\phi}$ analogously to $M_{\iota}$, $M_{\iota}^{\tau}$, and $M_{\iota}^{\phi}$, respectively.  
It follows from Proposition \ref{prop-n-even-independent-odd-ss-dependent}-(ii) that
\[
\det M_{\iota}^{s,\tau} = 0, \quad \text{for all } \iota = 1, \ldots, (m+1)n-1.
\]  
In particular, for $\iota = n$, this yields
\[
\det M_{n}^{s,\tau} = 0.
\]  
Substituting it into the identity
\[
\det M_{n}^{s} = \det M_{n}^{s,\tau} + \det M_{n}^{s,\phi}
\]
immediately gives
\[
\det M_{n}^{s} = \det M_{n}^{s,\phi}.
\]

Applying the same methodology as in the proof of part (ii) then produces the asserted expression for $\det M_{n}$.  

Finally, we verify that $\beta_s \neq 0$. Indeed, $\beta_s$ is the determinant of the submatrix obtained by removing the last row and the $n$-th column from $M^{s}$. Proposition \ref{prop-n-even-independent-odd-ss-dependent}-(ii) ensures that this submatrix is nonsingular, and hence $\beta_s \neq 0$.
    
\hfill$\Box$


\noindent
\textbf{Acknowledgements}
The project is partially supported by National Natural Science Foundation of China (Grant No. 12271038, 12526205), and the Open Project Program of Key Laboratory of Mathematics and Complex System (Grant No. K202503), Beijing Normal University.


\end{document}